\documentclass[12pt]{article}
\usepackage{amsmath,amsthm,amsfonts}

\newtheorem{thm}{Theorem}[section]
\newtheorem{lemma}[thm]{Lemma}
\newtheorem{prop}[thm]{Proposition}
\newtheorem{cor}[thm]{Corollary}
\newtheorem{defin}[thm]{Definition}
\newtheorem{rem}[thm]{Remark}

\newcommand{\prlabel}[1]{\label{#1}}
\newcommand{\R}{{\mathbb{R}}}

\newcommand{\Z}{{\mathbb{Z}}}
\newcommand{\N}{{\mathbb{N}}}
\newcommand{\C}{{\mathbb{C}}}

\newcommand{\La}{{\Lambda}}
\newcommand{\tLa}{{\widetilde{\Lambda}}}
\newcommand{\tP}{{\widetilde{\mathcal{P}}}}

\newcommand{\hJ}{{\hat{J}}}

\newcommand{\cA}{{\mathcal{A}}}
\newcommand{\cD}{{\mathcal{D}}}
\newcommand{\cE}{{\mathcal{E}}}
\newcommand{\cF}{{\mathcal{F}}}
\newcommand{\cG}{{\mathcal{G}}}
\newcommand{\cH}{{\mathcal{H}}}
\newcommand{\cI}{{\mathcal{I}}}
\newcommand{\cM}{{\mathcal{M}}}

\newcommand{\cP}{{\mathcal{P}}}
\newcommand{\cS}{{\mathcal{S}}}

\newcommand{\cV}{{\mathcal{V}}}
\newcommand{\cW}{{\mathcal{W}}}

\newcommand{\oQ}{{\overline{QH}_{ev} (M)}}

\newcommand{\tG}{{\widetilde{G}}}
\newcommand{\tmu}{{\tilde{\mu}}}
\newcommand{\tf}{{\tilde{f}}}

\newcommand{\tih}{{\tilde{h}}}

\newcommand{\Ca}{{\hbox{\it Cal}}}
\newcommand{\tCa}{{\widetilde{\hbox{\it Cal}}}}

\newcommand{\Ham}{{\hbox{\it Ham\,}}}
\newcommand{\spec}{{\text{spec}}}
\newcommand{\id}{{\text{{\bf 1}}}}
\def\square{{\vrule height6pt width6pt depth2pt}}
\begin{document}

\title{Calabi quasimorphism and quantum homology}

\author{\textsc Michael Entov
\and Leonid Polterovich\thanks{Supported by the Israel Science
Foundation}}

\date{September 2002}

\maketitle

\centerline{{\it Dedicated to the memory of Robert Brooks}}

 \begin{abstract}
\noindent
We prove that the group of
area-preserving  diffeomorphisms
of the 2-sphere admits a non-trivial
homogeneous quasimorphism to the real numbers
with the following property. Its value
on any diffeo\-mor\-phism supported in a
sufficiently small open subset of the sphere
equals to the Calabi invariant of the diffeo\-mor\-phism.
This result extends to more general symplectic manifolds:
If the symplectic manifold is monotone and its
quantum homology algebra is semi-simple we construct a
similar quasimorphism on the universal cover of the
group of Hamiltonian diffeomorphisms.
\end{abstract}

\vfill\eject \tableofcontents\vfill\eject

\section{Introduction and results} \prlabel{1}

A quasimorphism on a group $G$ is a function $r: G \to \R$ which
satisfies the homomorphism equation up to a bounded error: there
exists $R > 0$ such that $$|r(fg) -r(f) -r(g)| \leq R$$ for all
$f,g \in G$ (see \cite{Bav} for preliminaries on quasimorphisms).
A quasimorphism $r_h$ is called {\it homogeneous} if $r_h(g^m) = m
r_h(g)$ for all $g \in G$ and $m \in \Z$. Homogeneous
quasimorphisms are invariant under conjugations in $G$. Every
quasimorphism $r$ gives rise to a homogeneous one $$r_h(g) =
\lim_{m \to +\infty} \frac{r(g^m)}{m},$$ called the {\it
homogenization} of $r$. Starting from the classical work of
R.Brooks \cite{Brooks}, who explicitly constructed a non-trivial
quasimorphism on a free group, quasimorphisms are playing an
important role in the study of groups. In particular, they appear
in the bounded cohomology theory and in the geometry of the
commutator norm (see e.g. \cite{Bav} and Section \ref{1.3.A}
below).

In the present paper we focus on the cases when $G$ is either the
group of Hamiltonian diffeomorphisms of a symplectic manifold or
its universal cover. For a class of manifolds, which for instance
includes complex projective spaces, we give an explicit
construction of a non-trivial quasimorphism on $G$. The
construction involves some tools from "hard" symplectic topology,
in particular, Floer and quantum homology. Interestingly enough,
our quasimorphism is closely related to the classical Calabi
invariant.

\subsection{Extending the Calabi homomorphism} \prlabel{1.1}
Let $G = \Ham (M,\omega)$ be the group of Hamiltonian
diffeomorphisms of a closed connected symplectic manifold
$M^{2n}$ (see e.g.
\cite{McD-Sal-sympl-top},\cite{Pol-book} for preliminaries on $G$).
The group $G$ has a natural class of subgroups $G_U$ associated
to non-empty open subsets $U \subset M$, $U \neq M$. The subgroup $G_U$
consists of all elements $f \in G$ generated by a time-dependent
Hamiltonian
\begin{equation}
\prlabel{1.1.A}
F_t : M \to \R, \;t \in [0;1] \;\text{with}\;
\text{support}(F_t) \subset U\;
\text{for}\; \text{all}\; t.
\end{equation}
Consider the map $\Ca _U : G_U \to \R$ given by
\begin{equation}
\prlabel{1.1.B}
f \mapsto \int_0^1 dt \int_M F_t \omega^n.
\end{equation}
When the symplectic form $\omega$ is exact on $U$, this map is
well defined, meaning that it does not depend on the specific choice
of the Hamiltonian $F$ generating $f$. In fact, $\Ca _U$
is a homomorphism called {\it the Calabi homomorphism} \cite{Ban},
\cite{Cal}.
Note that $G_U \subset G_V$ for $U \subset V$, and in this case
$\Ca _U = \Ca _V$ on $G_U$.

In what follows we deal with the class $\cD$ of all non-empty
open subsets $U$
which can be {\it displaced} by a Hamiltonian
diffeomorphism:
\begin{equation}
\prlabel{1.1.C}
hU \cap \text{Closure}~(U)
= \emptyset \;\text{for}\;\text{some}\; h \in G.
\end{equation}
Put $$\cD_{ex} =
\{U \in \cD \;\big{|}\; \omega \; \text{is}\;
\text{exact}\;\text{on}\;U\}.$$

A celebrated result due to A.Banyaga \cite{Ban} states that the
group $G$ is simple and therefore does not admit a non-trivial
homomorphism to $\R$. In this paper we observe the following
phenomenon: for certain symplectic manifolds the family of
homomorphisms $\{\Ca _U: G_U \to \R\}_{U \in \cD_{ex}}$ extends to
a quasimorphism from $G$ to $\R$.

Given a Hamiltonian $H: M\times
S^1\to\R$ denote by $\psi_H\in G$ the Hamiltonian
symplectomorphism generated by $H$, i.e. time-1 map of the
Hamiltonian flow of $H$.

\begin{defin}
\label{def-cont-quasi}
{\rm
Suppose that a function $r : G\to \R$
satisfies the following condition: if a sequence
$\{ H_i\}$ of
smooth (possibly time-dependent) Hamiltonians $H_i: M\times
S^1\to \R$ converges $C^0{\hbox{\rm -uniformly}}$ to a smooth
function $H: M\times S^1\to\R$ then
\begin{equation}
\prlabel{A1.A}
r (\psi_{H_i}) \to r (\psi_H) \; \text{as} \; i \to \infty.
\end{equation}
In such a case the function $r$ will be called {\it continuous}.
}
\end{defin}

\begin{defin}
\label{def-Calabi-quasi}
{\rm
A quasimorphism on $G$ coinciding with the Calabi homomorphism
$\{\Ca _U: G_U \to \R\}$ on any ${U \in \cD_{ex}}$ will be called
{\it a Calabi quasimorphism}.
}
\end{defin}

\begin{thm}\prlabel{1.1.D}
Let $(M,\omega)$ be one of the following symplectic manifolds:
\begin{itemize}
\item{} the 2-sphere $S^2$ with an area form $\omega$;
\item{}$S^2 \times S^2$ with
the split symplectic form $\omega \oplus \omega$;
\item{} the
complex projective space $\C P^n$ endowed with the Fubini-Study
form.

\end{itemize}
Then there exists a continuous homogeneous Calabi
quasimorphism $\mu: G \to \R$.
\end{thm}

\bigskip
For all the cases except $\C P^n$, $n\geq 3$,
such a quasimorphism $\mu$ will be constructed in
Section~\ref{3} below.
For the case  $\C P^n$, $n\geq 3$, see
Section~\ref{sect-cp-n-proof}.

\begin{rem}
\label{rem-thm-G}
{\rm

\medskip
\noindent
1)  We do not know whether such a quasimorphism $\mu$ is unique.

\medskip
\noindent
2) In the case of $S^2$ one can show that
all continuous homogeneous Calabi quasimorphisms on $G$
coincide on elements of the group generated by auto\-nomous Hamiltonians:
given such a quasimorphism $\mu$ and an autonomous Hamiltonian $H:
S^2 \to\R$, one can explicitly compute $\mu (\psi_H)$ in terms of
combinatorics of the level sets of $H$ -- see Section~\ref{A}.

\medskip
\noindent
3) The {\it specific} homogeneous quasimorphism $\mu: G\to \R$
constructed in the proof of Theorem~\ref{1.1.D} is not only
continuous but also Lipschitz with respect to the famous Hofer
metric on $G$ -- see Section~\ref{3.7}.
}
\end{rem}

\bigskip
In fact, the "natural environment" in which one can look for a
Calabi quasimorphism is the universal cover $\tG$ of $G$ rather
then $G$ itself. Namely, for a non-empty open subset $U \subset
M$, $U \neq M$, consider a subgroup $\tG_U \subset \tG$
defined as follows.
An element of $\tG$ lies in $\tG_U$ if and only if it can be
represented by a Hamiltonian
flow $\{f_t\}_{t\in[0;1]}$, with $f_0 = \id$, generated by
a Hamiltonian $F_t$ satisfying condition (\ref{1.1.A}).
Formula
(\ref{1.1.B}) gives rise to a well defined homomorphism $$\tCa :
\tG _U \to \R.$$ We wish to extend the family of Calabi
homomorphisms $$\{\tCa _U\}_{U \in \cD}$$ to a
quasimorphism on $\tG$. If such an extension is possible the
resulting quasimorphism on $\tG$ is also called {\it a Calabi
quasimorphism}. The definition of a continuous function
on $\tG$ virtually repeats Definition~\ref{def-cont-quasi}.

The mere existence of a continuous
homogeneous Calabi quasimorphism on
$\tG$ can be shown
for a larger class of symplectic manifolds than the previous
theorem (see below). For the manifolds listed in Theorem
\ref{1.1.D} such a quasimorphism on $\tG$ actually descends
to $G$: for the cases other than $\C P^n$ ($n\geq 3$) this
is due to finiteness of the fundamental group $\pi_1 (G)$; in
the case $\C P^n$ ($n\geq 3$) when $\pi_1 (G)$ is
unknown, the proof relies on a delicate argument due to
P.Seidel and based on his work \cite{Se}.

Now we are going to formulate a result concerning the existence of
a homogeneous Calabi quasimorphism on $\tG$. It will hold for so called
spherically monotone symplectic manifolds. Recall that a closed
connected symplectic manifold $(M,\omega)$ is called {\it
spherically monotone} if there exists a real constant $\kappa > 0$ such
that $$(c_1 (M),A) = \kappa\cdot([\omega],A) \;
\text{for}\;\text{all}\; A \in \pi_2(M).$$ Here $c_1 (M)$ stands
for the first Chern class of the symplectic bundle $TM \to M$
equipped with an {\it $\omega$-compatible} almost complex
structure $J$ on $M$, where $\omega$-compatibility means that the form
$\omega (\cdot, J\cdot)$ is a Riemannian metric on $M$ (such an
almost complex structure is homotopically unique \cite{Gro-pshc}).

A crucial character of our story is the even-dimensional quantum
homology algebra $QH_{ev}(M)$ (see \cite{Liu},
\cite{McD-Sal-pshc}, \cite{Ru-Ti}, \cite{Ru-Ti-1}, \cite{Wi}) over
the field $k = \C[[s]$. Elements of $k$ are formal Laurent series
$\sum_{j \in \Z} z_j s^j$ where $z_j \in \C$, $s$ is a formal
variable and all $z_j$ vanish for large enough positive $j$. The
even-dimensional quantum homology is a commutative Frobenius
algebra with unity whose vector space structure is given by
$H_{ev}(M) \otimes_\C k$.\footnote{By $H_\ast (M)$ we always
denote the singular homology groups of $M$ with complex
coefficients, and $H_{ev} (M)$ stands for its even part.} The
product on $QH_{ev} (M)$ is a certain deformation of the
homological intersection product $$\cap: H_{ev}(M) \otimes
H_{ev}(M) \to H_{ev}(M).$$ Set $P=[\text{point}]\in H_0 (M)$. The
Frobenius structure on $QH_{ev} (M)$ is given by a non-degenerate
bilinear $k{\hbox{\rm -valued}}$ form $\Delta$ which associates to
a pair of quantum homology classes $a, b\in QH_{ev}(M)$ the
coefficient at $P$ in their quantum product $a\ast b \in H_{ev}(M)
\otimes_\C k$. We refer to \ref{2.3} below for brief preliminaries
on quantum homology and to \cite{McD-Sal-pshc} for a detailed
exposition.

Recall that a commutative algebra $Q$ over a field $k$ is called
{\it semi-simple} if it splits into a direct sum of fields as
follows: $Q = Q_1 \oplus...\oplus Q_d\;$, where
\begin{itemize}
\item{}each $Q_i \subset Q$
is a finite-dimensional linear subspace over $k$;
\item{} each $Q_i$ is a
field with respect to the induced ring structure;
\item{} The multiplication in $Q$ respects the splitting:
$$(a_1,...,a_d)\cdot(b_1,...,b_d) = (a_1 b_1,...,a_d b_d).$$
\end{itemize}

The semi-simplicity of a Frobenius algebra
$Q$
over a field
$k$
can be checked using a criterion due to L.Abrams
\cite{Abrams}
which says that a Frobenius algebra is semi-simple if and only if
its {\it Euler class}
is invertible. Recall that the Euler class
${\cal E}$
of a Frobenius algebra
$Q$
is defined as
\[
{\cal E}= \sum_i e_i e_i^\sharp,
\]
where
$\{ e_i\}$
is a basis of
$Q$
over
$k$
and
$\{ e_i^\sharp \}$
is the dual basis with respect to the non-degenerate bilinear form
on
$Q$
defining the Frobenius structure.
The Euler class does not depend on the choice of the basis
$\{ e_i\}$.

\begin{thm}
\prlabel{1.1.E}
\label{thm-11E}
Let $(M,\omega)$ be a closed connected spherically mo\-no\-tone
symplectic manifold. Suppose that the quantum homo\-logy
algebra
$QH_{ev}(M)$ is semi-simple.
Then there exists a continuous
homogeneous Calabi quasi\-mor\-phism $\tilde{\mu}: \tG \to \R$.
\end{thm}

\bigskip
\begin{rem}
\label{rem-thm-tG}
{\rm
\medskip
\noindent
1) We do not know whether such a quasimorphism $\tilde{\mu}$ is
unique.

\medskip
\noindent 2) The quasimorphism $\tilde{\mu}$ can be calculated
on the subgroup $\pi_1 (G)\subset \tG$ in terms of the
Seidel action of $\pi_1(G)$ on the quantum homology of $M$ (see
Section~\ref{sect-c-on-loops}).

\medskip
\noindent
3) The Hofer metric on $G$ can be lifted to
a (bi-invariant) pseudo-metric on $\tG$.
The quasimorphism $\tilde{\mu}$ is Lipschitz with respect to this
pseudo-metric (see Section~\ref{3.7}).
}
\end{rem}

\bigskip
Examples of symplectic manifolds $M$ with semi-simple quantum homology
algebra $QH_{ev} (M)$ include, in particular,
$S^2$, $S^2\times S^2$,
$\C P^n$, $\C P^2$ blown up at
one point and complex Grassmannians
with the usual monotone
symplectic structures.
To get the semi-simplicity
of $QH_{ev} (M)$ in these cases one can use the
known explicit descriptions of the multiplicative structure of
$QH_\ast (M)$ to check that the Euler class is invertible so that
the Abrams criterion can be applied. For more details on the first
four examples see
Section~\ref{sec-qh-examples}. In the case of a complex Grassmannian the
structure of the quantum homology algebra
is described in
\cite{Bertr},
\cite{Sieb-Ti},
\cite{Wit-grasm}.
In such a case the Euler class is an integral multiple of
$P = [\text{point}]$ (see
\cite{Abrams},\cite{Bertram-Euler})
which is invertible according to a computation based on
\cite{Bertr-CiF-Fult} and due to A.Postnikov
(see
\cite{Entov}, also see \cite{Post-new}).

It is known that when the class of the symplectic form $[\omega]$
vanishes on $\pi_2(M)$ then the product structure
on $QH_{ev}(M) = H_{ev}(M) \otimes k$ is given by the ordinary
intersection product
$\cap$, and so is never semi-simple. Thus our result
does not apply to
those spherically monotone symplectic manifolds $(M, \omega)$ where
$[\omega]$ vanishes identically on
$\pi_2 (M)$.

\subsection{Applications and discussion} \prlabel{1.3}

\subsubsection{The commutator norm}\prlabel{1.3.A}
Let $G$ be a group and $[G, G]$ be its commutator subgroup. Every
element $h \in [G, G]$ can be written as a product of simple
commutators $fgf^{-1}g^{-1},\; f,g \in G$. The {\it commutator
norm} $|| h ||$ is by definition the minimal number of simple
commutators needed in order to represent $h$. It is known (see
e.g. \cite{Ba-Ghys}, \cite{Bav}) that for a homogeneous
quasimorphism $\mu : G \to \R$ one has \footnote {Here and below
$\hbox{\it const}$ stands for a positive constant.} $$|| h || \geq
\hbox{\it const}\, (\mu)\cdot \mu(h),\;\; h \in [G,G].$$ In
particular, existence of a homogeneous quasimorphism which does
not vanish on $[G,G]$ implies that the diameter of the group
$[G,G]$ with respect to the commutator norm is infinite. Recall
that $G$ is called {\it perfect} if $G = [G, G]$.

A.Banyaga \cite{Ban} proved that the group $G = \Ham(M,\omega)$ and
its universal cover
$\tG$ are perfect for every closed symplectic manifold
$(M,\omega)$.
As an immediate consequence of our results we get the following

\begin{cor}
\prlabel{1.3.B} Let $(M,\omega)$ be a closed connected spherically
mo\-no\-tone symplectic manifold, $G=\Ham (M,\omega)$. Suppose
that the quantum homology algebra $QH_{ev}(M)$ is semi-simple. Let
$U \in \cD$ be a displaceable open subset. Then $$||\tf || \geq
\hbox{\it const} \cdot |\tCa\, (\tf)| $$ for every $\tf \in \tG
_U$. If in addition the fundamental group $\pi_1 (G)$ is finite
then $$|| f || \geq \hbox{\it const} \cdot |\Ca\, (f)| $$ for
every $f \in G_U$ provided $U \in \cD_{ex}$.
\end{cor}

\bigskip
The second part of the corollary follows from
Proposition~\ref{3.6.A} below which says that if
$\pi_1 (G)$ is finite then  the Calabi quasimorphism $\tmu$ on $\tG$
descends to a homogeneous Calabi quasimorphism $\mu$ on $G$.
Corollary~\ref{1.3.B} generalizes a result obtained in a recent work
\cite{Entov} which served as the starting point for the present
research.

\subsubsection{Quantitative fragmentation lemma}\prlabel{1.3.D}

 Let $\{U_1,...,U_m\}$ be an open
covering of a closed connected symplectic
manifold $(M,\omega)$. Banyaga's fragmentation lemma states
that any element $f \in G$ can be written as a product of
diffeomorphisms $g_i$ as follows.
Each $g_i$ lies in $G_{U_j}$ for some $j\in \{1;...;m\}$
and moreover it is contained in the kernel of the Calabi
homomorphism. Denote by $l (f)$ the minimal number of $g_i$'s
needed in order to represent $f$.

Suppose now that $(M,\omega)$ is one of the manifolds
listed in Theorem \ref{1.1.D}, and all the sets $U_j$
lie in $\cD_{ex}$. The following is an immediate consequence of
Theorem~\ref{1.1.D}.

\begin{cor}
\prlabel{1.3.E} $$l (f) \geq \hbox{\it const} \cdot |\Ca\, (f)| $$
for every $f \in G_U$ provided $U \in \cD_{ex}$.
\end{cor}

\subsubsection{Asymptotic growth of one-parametric
subgroups}\prlabel{1.3.BB}

We recall a few known definitions.  Denote by $\cF$ the space of all smooth
Hamiltonian functions
$F:M \times S^1 \to \R$ which satisfy the following normalization
condition:
$ \int_M F_t\, \omega^n = 0$
for all $t \in S^1$, where
$F_t = F(\cdot, t)$.
Introduce the $C^0$-norm on $\cF$ by
\begin{equation}
\label{conorm}
\| F\|_{C^0} = \max_M  F - \min_M F .
\end{equation}

A distance between the identity $\id$ and an
element $f$ of the group $G$ is defined \cite{Ho} as
\[
\rho (\id, f) = \inf_F \int_{S^1} \| F_t\|_{C^0}\;
dt,
\]
where the infimum is taken over all time-dependent Hamiltonians
$F\in\cF$ generating $f$.
The distance function $\rho$ gives rise to
a bi-invariant non-degenerate metric on $G$ \cite{Ho},
\cite{Lal-McD-metr}, \cite{Pol-metr}, called the {\it Hofer metric}.

A {\it time-independent} Hamiltonian $F\in \cF$
generates a one-parametric subgroup
$\{ \psi^t_F \}$ of $G$ so that $\psi_F = \psi^1_F$. The {\it asymptotic growth}
of the subgroup $\{ \psi^t_F \}$ is defined \cite{Pol-book} as
\[
\zeta (F) = \lim_{t\to +\infty} \frac{\rho
(\id, \psi_F^t)}{t \| F\|_{C^0}}.
\]
Such a limit always
exists and belongs to $[0,1]$.

\begin{cor}
\prlabel{1.3.EE}
Let $M=S^2$.
Then for a generic $F$
\[
\zeta (F) > 0.
\]
\end{cor}

\bigskip
The proof can be found in Section~\ref{sect-pf-cor-asympt-length}. It relies
on the
estimate
\begin{equation}
\label{eqn-zeta-mu-0}
\zeta (F) \geq \frac{|\mu (\psi_F) |}{\| F\|_{C^0}},
\end{equation}
which holds for the {\it specific} quasimorphism $\mu$ constructed
in the proof of Theorem~\ref{1.1.D} (see
Section~\ref{sect-pf-cor-asympt-length}), and on the explicit
computation of the value of $\mu: G\to\R$ on $\psi_F^m$,
$m=1,2,\ldots$, in the case $M=S^2$ made in Section~\ref{A4}.

\begin{rem}
\label{rem-asympt-growth-S-2}
{\rm

\bigskip
\noindent
1) In fact, the inequality $\zeta (F) > 0$ for a
generic Hamiltonian $F$ is valid on all closed symplectic
surfaces. The case of the 2-torus is settled in
\cite{Pol-book}, Section 8.4. The argument given in
\cite{Pol-book}
actually works for any closed symplectic
surface of positive genus.

\bigskip
\noindent
2) For symplectic manifolds listed in
Theorem~\ref{1.1.D}
inequality (\ref{eqn-zeta-mu-0}) immediately produces
examples of 1-parametric subgroups of $G$ with positive
$\zeta (F)$ -- for example, take an autonomous Hamiltonian $F$ supported
in a sufficiently small ball and such that the Calabi invariant of $\psi_F$
is non-zero. This shows, in particular, that for these manifolds
the group $G$ has infinite diameter with respect to the Hofer metric.
}
\end{rem}

\subsubsection{Other quasimorphisms?} \prlabel{1.3.C}
J.Barge and E.Ghys \cite{Ba-Ghys} constructed a quasimorphism
of a different nature on the group of compactly supported
symplectomorphisms of a standard symplectic ball.
The Barge-Ghys quasimorphism is closely related to the Maslov
class in symplectic geometry.
This construction was later
generalized in \cite{Entov} to
closed symplectic manifolds $(M, \omega)$ with $c_1 (M) =0$
(e.g. tori and K3-surfaces). As a result one gets
a homogeneous quasimorphism on
$\widetilde{\textrm Symp}_0\, (M, \omega)$
which does not vanish on
$\widetilde{\textrm Ham}\, (M, \omega)$ \cite{Entov}.
Here
$\widetilde{\textrm Symp}_0\, (M, \omega)$
is the universal cover of the identity component of the group of
symplectomorphisms of $(M, \omega)$.

Interestingly enough, this class of manifolds
is disjoint from the one
considered in the present paper.
We believe however that the class of manifolds
admitting a Calabi quasimorphism can be enlarged, namely
the spherical monotonicity condition can be removed.
Such a generalization
should go
\break
along the same lines though the technicalities
will become more complicated.
On the
other hand, semi-simplicity of the quantum homology algebra
seems to be a crucial assumption.

No other quasimorphism of $G$ and $\tG$ is known at the moment.
The simplest symplectic manifolds for which no information on
quasimorphisms and the commutator norm is available at all
are closed
oriented surfaces of
genus $\geq 2$. Any progress in this direction would
be very interesting.

\subsubsection{Continuum of Calabi quasimorphisms on an open
surface} \prlabel{1.2.5}

The notion of Calabi quasimorphism can be extended in a
straightforward way to {\it open}  symplectic manifolds. It turns
out that even for very simple manifolds Calabi quasimorphisms can
form an infinite-dimensional affine space. For an open connected
symplectic manifold $M$ denote by $G_M$ the group of all
Hamiltonian diffeomorphisms of $M$ generated by Hamiltonians with
compact support in $M \times S^1$. It is known that formula
(\ref{1.1.B}) gives rise to the well defined Calabi homomorphism $
\Ca_M : G_M \to \R$. Up to a multiple, this is the only
homomorphism $G_M \to \R$.

\begin{thm}
\prlabel{1.2.5.A} Suppose that either $M \subset \R^2$ is an open
disk of finite area, or $M \subset T^*S^1$ is an open annulus of
finite area. There exists a family $\mu_{\epsilon},\;\epsilon \in
\R$, of continuous homogeneous Calabi quasimorphisms on $G_M$ with
the following properties:
\begin{itemize}
\item{} Given a finite subset $I \subset \R$, the quasimorphisms
$\mu_{\epsilon},\; \epsilon \in I$, are linearly independent over
$\R$. In particular, the 2-nd bounded cohomology of $G_M$ is an
infinite-dimensional space over $\R$;
\item{} Moreover, if $M$ is an annulus, the quasimorphisms
can be chosen so that every $\mu_{\epsilon}$ coincides with the
Calabi homomorphism $\Ca_M$ on the subgroup $G_U$, where $U$ is
the interior of any (not necessarily displaceable!) embedded
closed disk in $M$.
\end{itemize}
\end{thm}

\medskip
\noindent The proof is given in Section \ref{bcoh} below. It seems
likely that analogous results hold true for any other open surface
of genus $0$ with finite area. It would be also interesting to
find a generalization to higher dimensions, for instance to the
cases when $M$ is either the standard symplectic open ball, or the
open unit coball bundle of the flat $n$-dimensional torus. Let us
mention also that J.-M. Gambaudo  \cite{Gam} suggested a different
approach which could lead to an infinite sequence of non-trivial
quasimorphisms on $G_M$ in the case when $M$ is a 2-dimensional
disk.

\section{Symplectic preliminaries}\prlabel{2}

\subsection{Starting notations}\prlabel{2.1}
Let $(M^{2n},\omega)$ be a closed
spherically monotone symplectic manifold.
Consider the group
$$\bar \pi_2 (M) = \pi_2(M) / \sim,$$
where by definition $A \sim B$ iff $([\omega], A) = ([\omega], B)$.
Clearly, both $[\omega]$ and $c_1 (M)$ descend to homomorphisms of
$\bar \pi_2 (M)$.
In view of the comment at the very end of Section~\ref{1.1}, we will
always assume that $[\omega]$ does not vanish on $\bar \pi_2(M)$.
In particular,
the group $\bar \pi_2 (M)$ is the infinite cyclic group,
and it has a generator $S$
so that
$\Omega :=([\omega], S) > 0$. Set $N := (c_1(M), S) >0$.
As above $\tG$ stands for the universal cover of $\Ham (M,\omega)$.

\subsection{Gromov-Witten invariants}\prlabel{2.2}
The Gromov-Witten invariant $GW_j,\; j \in \N ,$ is a 3-linear
(over $\C$) $\C$-valued form on $H_*(M)$.
It is defined along the following lines (see \cite{McD-Sal-pshc},
\cite{Ru-Ti}, \cite{Ru-Ti-1} for the precise
definition). First of all, $GW_j (A,B,C) = 0$ unless
$$\deg A + \deg B + \deg C = 4n - 2Nj.$$
If the equality above holds, we assume without loss of generality
that the homology classes $A,B$ and $C$
are represented by smooth submanifolds $\widehat{A}, \widehat{B}$ and $\widehat{C}$ respectively.
Take an $\omega$-compatible almost complex structure $J$ on $M$.
Consider the following elliptic problem:

\medskip
\noindent
Find all $J$-holomorphic maps $\C P^1 \to (M,J)$ which
represent the class $jS \in \bar \pi_2 (M)$ and which send the
points $0,1,\infty \in \bar{\C} \cup \{\infty\} = \C P^1$
to $\widehat{A},\widehat{B}$ and $\widehat{C}$ respectively.

\medskip
When the almost complex structure $J$ and the submanifolds $\widehat{A},\widehat{B},\widehat{C}$
are chosen in a generic way, the set of solutions of this problem is
finite. The
number $GW_j(A,B,C)$ is defined as the number of the solutions
counted with an appropriate sign. It is useful to have in mind that
if $J$ is a genuine
complex structure and $\widehat{A}$,$\widehat{B}$,$\widehat{C}$ are generic complex
submanifolds, the sign in question is positive.

\subsection{Quantum homology algebra} \prlabel{2.3}
\label{sec-qh-examples}
As a vector space over $\C$ the quantum homology $QH_\ast (M)$ is
isomorphic to $H_*(M) \otimes_{\C} k$, where $k$ stands for
the field
$\C[[s]$ which appeared in Section~\ref{1.1}.
The quantum multiplication $a*b$, $a,b\in QH_\ast (M)$, is defined as follows.
For $A,B \in H_*(M)$ and $j \in \N$
define $(A*B)_j \in H_*(M)$ as the unique class which
satisfies $$(A*B)_j \circ C = GW_j(A,B,C)$$
for all $C \in H_*(M)$. Here $\circ$ stands for the ordinary
intersection index in homology.
Now for any
$A, B\in H_\ast (M)$
set
$$A*B = A \cap B + \sum_{j \in \N} (A*B)_j\; s^{-j} \in QH_*(M).$$
By $k{\hbox{\rm -linearity}}$ extend the quantum product
to the whole $QH_\ast (M)$. As a result one
gets a correctly defined skew-commutative associative product
operation on $QH_\ast (M)$
which is a deformation of the classical $\cap$-product in singular
homology
\cite{Liu}, \cite{McD-Sal-pshc}, \cite{Ru-Ti}, \cite{Ru-Ti-1},
\cite{Wi}.

The field $k$ has a ring grading defined by the condition that the grade
degree of $s$ equals $2N$.
Such a grading on $k$ together with
the usual grading on
$H_\ast (M)$
define a grading on
$QH_\ast (M) = H_\ast (M)\otimes_{\C} k$.
If $a, b\in QH_\ast (M)$ have graded degrees $deg\, (a)$, $deg\, (b)$
then
$deg\, (a\ast b) = deg\, (a) + deg\, (b) - 2n$.

The fundamental class $[M]$ is the unity with respect to
the quantum
multiplication. If $A\in H_\ast (M)$, $\varsigma\in k$, we
will denote the elements
$A\otimes 1, [M]\otimes \varsigma$ of $QH_\ast (M) = H_\ast (M)\otimes_{\C} k$
respectively by $A$ and
$\varsigma$.
The even part $QH_{ev} (M) := H_{ev}(M) \otimes_{\C} k$ is a commutative
subalgebra of $QH_\ast (M)$.

The algebra $QH_{ev}(M)$ is a Frobenius algebra over $k$. Consider
the pairing $\Delta: QH_{ev}(M) \times QH_{ev}(M)  \to k$,
$$\Delta {\big (}\sum A_j s^j, \sum B_l s^l{\big )} = \sum (A_j
\circ B_l) \cdot s^{j+l}.$$ In fact $\Delta$ associates to a pair
of quantum homology classes $a, b\in QH_{ev}(M)$ the coefficient
at $P = [\text{point}]$ in their quantum product $a\ast b \in
H_{ev}(M) \otimes_\C k$.  The pairing $\Delta$ defines a Frobenius
algebra structure , which means that $\Delta$ is non-degenerate
and $$\Delta(a,b) = \Delta(a*b,[M])\;\;
\text{for}\;\text{all}\;a,b \in QH_{ev}(M).$$ Let $\tau: k \to \C$
be the map sending $\sum z_j s^j$ to $z_0$. Define a $\C$-valued
pairing
\begin{equation}
\label{3.2.A}
\Pi (a,b) = \tau\Delta(a,b) = \tau\Delta(a*b,[M])
\end{equation}
on $QH_{ev}(M)$. It will play an important role below.

\medskip
\noindent
\subsubsection{Example: $S^2$}\prlabel{2.3.A}
Let $M$ be the 2-sphere $S^2$. Set $P = [\text{point}]$.
Note that $GW_1 (P,P,P) =1$, and moreover this is the only non-vanishing
Gromov-Witten invariant. Hence $P*P = s^{-1}$ and
$QH_\ast (M)$ is a field:
\[
QH_\ast (M) = \frac{k [P]}{\{ P^2 =s^{-1}\} }.
\]
Thus
$QH_\ast (M) = QH_{ev} (M)$
is semi-simple.

\subsubsection{Example: $\C P^n$}\prlabel{2.3.D}

The previous example
can be generalized.
Let
$M = {\C} P^n$
be equipped with the standard Fubini-Study symplectic form.
This is a spherically mono\-tone symplectic manifold.
Let
$A\in H_{2n-2} (M)$
be the projective hyperplane class.
Then
\[
QH_\ast (M) = \frac{k [A]}{ \{ A^{n+1} =  s^{-1}\} },
\]
(see
\cite{Ru-Ti},
\cite{Ru-Ti-1},
\cite{Wi}).
One immediately sees that
$QH_\ast (M) = QH_{ev} (M)$
is a field over
$k = \C [[s]$
and therefore it is a
semi-simple algebra.

\subsubsection{Example: $S^2\times S^2$}\prlabel{2.3.C}

Let
$M = S^2\times S^2$
be equipped with the split symplectic form
$\omega\oplus\omega$,
where
$\omega$
is an area form on
$S^2$.
This is a spherically monotone symplectic manifold.
Let $A, B\in H_2 (M)$
be the homology classes of
$S^2\times pt$
and
$pt\times S^2$.
The elements
$P = [\text{pt}], A, B, [M]$
form a basis of
$QH_\ast (M)$
over
$k$
and the multiplicative relations are completely defined by the identities:
\[
A\ast B = P,\ \  A^2 = B^2 = s^{-1}.
\]
Thus the Euler class is
\[
{\cal E} = 2 P\ast [M] + 2 A\ast B = 4 P.
\]
The classes $A, B$
are invertible and so are
$P$
and
${\cal E}$.
Therefore, according to the Abrams criterion, the algebra
$QH_\ast (M)  = QH_{ev} (M)$
is semi-simple.
Note that $QH_{ev} (M)$ is not a field since it contains divisors
of zero: $(A-B)\ast (A+B) = A^2 - B^2 = 0$.

\subsubsection{Example: $\C P^2$ blown up at one point}\prlabel{2.3.E}

Let
$M$
be the complex blow up of
$\C P^2$
at one point equipped with a monotone symplectic form
(see
\cite{McD-geom-var},\cite{Pol-masl}).
Its quantum homology algebra is described in
\cite{McD-geom-var}
as follows. Let
$P = [\text{point}]$.
Denote by
$A$
the exceptional divisor and set
$B = [\C P^1] -A$.
Together with the fundamental class
$[M]$
the classes
$P, A, B$
generate
$QH_\ast (M)$
as a vector space over
$k$.
The multiplicative relations are as follows
(recall that $[M]$ is the unity element in the quantum homology
algebra):
\[
\begin{array}{lll}
P\ast P  =  (A+B) s^{-3} &\ \ \ \ &A \ast P  =  B s^{-2}\\
P\ast B  =   s^{-3} & \ \ \ \ &A\ast A  =  -P + A s^{-1} +  s^{-2}\\
A\ast B  =  P - A s^{-1}  &\ \ \ \ & B\ast B  = A s^{-1}.\\
\end{array}
\]
The Euler class ${\cal E}$ is easily computable:
\[
{\cal E} =
P\ast [M] + A\ast B + B\ast (A +B) +[M]\ast P = 4P -  A s^{-1}.
\]
One can check that ${\cal E}$ is invertible:
\[
{\cal E}^{-1} = \frac{1}{283} (-12 P s^4 +9 A s^3 +73 B s^3 +
16 s^2).
\]
Therefore, according to the Abrams criterion,
$QH_\ast (M) = QH_{ev} (M)$
is a semi-simple algebra.

\subsection{The action functional}\prlabel{2.4}
Let $\La$ be the space of all smooth contractible loops $x: S^1 =
\R/\Z \to M$. Consider a covering $\tLa$ of $\La$ whose elements
are equivalence classes of pairs $(x,u)$, where $x \in \La$,  $u$
is a disk spanning $x$ in $M$ and the equivalence relation is
defined as follows: $(x_1,u_1) \sim (x_2,u_2)$ iff $x_1 = x_2$ and
the 2-sphere $u_1 \# (-u_2)$ vanishes in $\bar \pi_2 (M)$. The
equivalence class of a pair $(x, u)$ will be denoted by $[x,u]$.
The group of deck transformations of this covering can be
naturally identified with $\bar \pi_2 (M)$. The generator $S$ acts
by the transformation $s: \tLa \to \tLa$ so that
\begin{equation}
\label{act1}
s ([x,u]) = [x,u \# (-S)].
\end{equation}

Recall that by $\cF$ we denote the space of all smooth Hamiltonian functions
$F:M \times S^1 \to \R$ which satisfy the following normalization
condition:
$ \int_M F (\cdot , t)\, \omega^n = 0$
for any $t \in S^1$. For $F \in \cF$ define {\it the action functional}
$$\cA_F ([x,u]) := \int_{S^1} F(x(t),t) dt - \int_u \omega$$
on $\tLa$. Note that
\begin{equation}
\label{act2}
\cA_F(sy) =\cA_F(y) + \Omega
\end{equation}
for all $y \in \tLa$.

Let $\cP_F \subset \La$ be the set
of all contractible 1-periodic orbits of the
Hamiltonian flow
generated by $F$. Its full lift $\tP_F$ to $\tLa$ coincides
with the set of critical
points of $\cA_F$. We define the {\it action spectrum}
$\spec (F)$ as the set of
critical values of $\cA_F$. This is a closed nowhere dense
$\Omega \Z$-invariant subset of $\R$ \cite{Oh-act}, \cite{Sch-1}.

\subsection{Filtered Floer homology}\prlabel{2.5}
For a generic Hamiltonian $F \in \cF$ and $\alpha \in (\R \setminus \spec(F)) \cup \{+\infty \}$
define a complex vector space  $C_{\alpha} (F)$ as the set of all formal
sums
$$\sum_{y \in \tP_F} z_y y,$$
where $z_y \in \C$, $\cA_F(y) < \alpha$, which satisfy the following
finiteness condition:
$$\# \{y\;|\; z_y \neq 0 \; \text{and} \; \cA_F(y) > \delta \}
< \infty$$
for every $\delta\in\R$.
Formula (\ref{act1}) defines a
structure of the vector space over $k$ on
$C_{\infty}(F)$.

Given a loop $\{J_t\},\; t\in S^1$, of $\omega$-compatible
almost complex structures, define a Riemannian metric on $\La$
by
$$(\xi_1,\xi_2) = \int_0^1 \omega (\xi_1 (t), J_t\xi_2(t))dt, $$
where $\xi_1,\xi_2 \in T\La$. Lift this metric to $\tLa$ and consider
the negative gradient flow of the action functional $\cA_F$.
For a generic choice of the Hamiltonian $F$ and the loop
$\{J_t\}$ the count of isolated
gradient trajectories connecting critical points of $\cA_F$
gives rise in the standard way \cite{Floer}, \cite{Ho-Sa}
to a Morse-type differential
\begin{equation}
\label{eqn-D-def}
d: C_{\infty}(F) \to C_{\infty}(F),\; d^2 =0.
\end{equation}
The differential $d$ is $k$-linear. Moreover it preserves $\C$-subspaces
$C_{\alpha}(F) \subset C_{\infty}(F)$ for all $\alpha \in \R$.

The complexes $(C_{\alpha}(F), d)$
have a natural grading according to
the Conley-Zehnder index ${\textit ind}: \tP_F: \to \Z$
(see \cite{Co-Ze}) which satisfies ${\textit ind}\, (sy) = {\textit ind}\, (y) + 2N$ for every
$y \in \tP_F$.
Note that different authors
use slightly different versions of the Conley-Zehnder
index. In order to fix our convention consider the case of
a sufficiently
$C^2{\hbox{\rm -small}}$
autonomous Morse Hamiltonian
$F$.
Then the Conley-Zehnder index
${\textit ind}\, (y)$
of an element
$y\in \tP_F$,
represented by a pair $(x, u)$
consisting of a critical point
$x$
of
$F$
and the trivial disk
$u$,
is equal to the Morse
index
of
$x$.
In what follows we are interested in
the {\it even} part of the homology of these complexes.

Notice that in spite of the involvement of the almost complex
structures $\{J_t\},\; t\in S^1$, in the definition of the complex
$(C_{\infty}(F), d)$,
 different choices of $\{J_t\}$
lead to complexes whose homologies are related by natural
isomorphisms {\it preserving the filtration} as long as the
Hamiltonian $F$ is fixed \cite{Floer}, \cite{Ho-Sa},
\cite{Oh-act}, \cite{Oh-new}. Because of this we will suppress the dependence
on $\{J_t\}$, and define
$$V_{\alpha}(F) = H_{ev}(C_{\alpha}(F),d) \;\text{and}\;
V^{\alpha}(F) = H_{ev} (C_{\infty}(F)/C_{\alpha}(F), d).$$

Meanwhile these homology groups have
been defined for generic Hamiltonians $F$
only.
Using an appropriate continuation  procedure
one can extend the definition to all $F \in \cF$. Namely, let
$\alpha\notin \spec(F)$.
Pick any generic Hamiltonians
$F^\prime, F^{{\prime\prime}}$
such that the spaces
$V_\alpha ({F^\prime})$
and
$V_\alpha ( F^{{\prime\prime}})$
are defined.
Then if
$F^\prime, F^{{\prime\prime}}$
are sufficiently $C^\infty{\textrm{-close}}$ to $F$ we have that
$\alpha\notin \spec(F^\prime)$,
$\alpha\notin \spec(F^{{\prime\prime}})$,
and the spaces
$V_\alpha ({F^\prime})$
and
$V_\alpha (F^{{\prime\prime}})$
are canonically isomorphic. Thus we can set $V_\alpha (F)$
as $V_\alpha ({F^\prime})$ for any $F^\prime$ sufficiently close to
$F$, and this definition is correct.
Similarly one can define $V^{\alpha}(F)$
for arbitrary $F$ and $\alpha \notin \spec (F)$.

Suppose now that two Hamiltonian functions $F,F' \in \cF$ generate
the same element $f \in \tG$. In this case $\spec (F) = \spec
(F')$. This was proved in \cite{Sch-1} (see Lemma 3.3 there) in
the case when $\omega$ vanishes on $\bar{\pi}_2$ and the proof
readily extends to the general case (see e.g. \cite{Oh-new}). The
resulting set will be denoted simply by $\spec (f)$. Moreover the
vector spaces $V_{\alpha}(F)$ and $V^{\alpha}(F)$ can be
canonically identified, respectively, with $ V_{\alpha}(F')$ and
$V^{\alpha}(F')$. Therefore we shall denote them, respectively, by
$V_{\alpha}(f)$ and $V^{\alpha}(f)$, where $\alpha \in (\R
\setminus \spec(f) )\cup {+\infty}$ (the dependence on $\{ J_t\}$
is suppressed for the same reasons as above). These spaces are
called {\it filtered Floer homology} of an element $f \in \tG$.

\subsection{Algebraic data and spectral invariants}\prlabel{2.6}

Filtered Floer homology come with additional algebraic data
which we are going to list now.

\subsubsection{Identification with quantum homology}\prlabel{2.6.A}

All spaces $V_{\infty} (f)$ are canonically identified (as vector
spaces over $k$) with $QH_{ev}(M)$ \cite{PSS}. The different
choices of a Hamiltonian $F$ generating $f$ and of almost complex
structures $\{ J_t\}$ give rise to different Floer complexes whose
homologies are related by natural isomorphisms preserving the
filtration, thus leading to a well-defined space $V_{\infty} (f)$.
For each such Floer complex its homology can be canonically
identified with $QH_{ev}(M)$ (as in \cite{PSS}). These
identifications agree for different Floer complexes and therefore
lead to a well-defined identification of $V_{\infty} (f)$ with
$QH_{ev}(M)$ which  preserves the grading \cite{PSS}. Recall that
$QH_\ast (M) = H_\ast (M)\otimes_\C k$ carries a grading, while
$V_{\infty} (f)$ is graded by the Conley-Zehnder index ${\textit
ind}$.

\subsubsection{Natural inclusions}\prlabel{2.6.B}
For any $\beta < \alpha \leq +\infty$ the natural inclusion
$C_{\beta}(F) \to C_{\alpha}(F)$ of Floer complexes leads to a
homomorphism
$$i_{\alpha\beta}: V_{\beta}(f) \to V_{\alpha}(f)$$
between their homology. Moreover $i_{\alpha\beta}i_{\beta\gamma} =
i_{\alpha\gamma}$ for any $\gamma < \beta < \alpha$. We abbreviate
$i_{\alpha}$ for $i_{\alpha\infty}$. In the case we wish to
emphasize the dependence of $i_{\alpha}$ on the element $f$ we
will write $i_{\alpha}\{f\}$.

\subsubsection{Spectral invariants: finiteness and continuity}\prlabel{2.6.C}

Following the works of C.Viterbo
\cite{Vit}, Y.-G. Oh \cite{Oh1}, \cite{Oh2}, \cite{Oh-act}, \cite{Oh-new} and
M.Schwarz \cite{Sch-1}, \cite{Sch-talks},
we give the following definition.
Given $f\in\tG$ and $a \in QH_{ev}(M)=V_{\infty}(f),\; a \neq 0$, set
$$c(a,f) = \inf\, \{\,\alpha\;|\; a \in \text{Image}\; i_{\alpha}\;\}.$$
Then $ -\infty < c(a,f) < +\infty$
and for a given $a$ the function
$c(a,f)$ is continuous with respect to the $C^{\infty}$-topology
on $\tG$.
In the case when
$a$
is a singular homology class this has been proved in
\cite{Oh-act}. The proof for the general case can be found in
\cite{Oh-new}.
In fact, in the case of a spherically monotone symplectic manifold,
when we consider an arbitrary quantum homology class,
the only detail that should be added to the proof in
\cite{Oh-act}
is the following one.

A generic element $f\in \tG$ can be defined by means of a
Hamiltonian flow (generated by a Hamiltonian $F$) that has only a
finite number of 1-periodic trajectories.
Thus, since
$M$
is spherically monotone, there exist
constants $R_1, R_2 > 0$ such that for any
$y\in \tP_F$
\begin{equation}
\label{eqn-action-index}
R_1\cdot  {\textit ind}\, (y) - R_2 \leq   \cA_F (y) \leq R_1 \cdot
{\textit ind}\, (y) + R_2.
\end{equation}
Now write  $a$
as a sum $a = \sum_m a^{(m)}$ of its
homogeneous graded components, with each
$a^{(m)}$
having the grading $m$. Consider the set $\cI := \{m\;|\; a^{(m)} \neq 0\}$.
The definition of the field $k$ implies that
$m_0 := \max \cI  < \infty$.
We are going to use the grading-preserving identification
of $QH_{ev} (M)$ with $V_\infty (f)$ (see
Section~\ref{2.6.A}).

First, it follows that
\begin{equation}
\label{eqn-c-max}
c (a, f) = \max_{m\in \cI} c (a^{(m)}, f ).
\end{equation}
Further, the classes $a$ and $a^{(m)}$
can be viewed as Floer homology classes from
$V_{\infty} (f)$.
Let
$C^{(m)} = \sum_y z_y y$,
$z_y\in \C$, $y\in \tP_F$,
be a Floer chain representing $a^{(m)}$.
Then
${\textit ind}\, (y) =m$ for any $y$
entering
$C^{(m)}$.
Using
(\ref{eqn-action-index})
we see that for any such $y$
\[
R_1 m - R_2 \leq  \cA_F (y) \leq R_1 m + R_2.
\]
Therefore, since the Floer chain $C^{(m)}$ representing $a^{(m)}$ was
chosen arbitrarily,
\[
R_1 m - R_2 \leq  c(a^{(m)}, f) \leq
 R_1 m + R_2 \; \text{for}\;\text{every}\; m \in \cI.
\]
Combining it with
(\ref{eqn-c-max}) we see that
\[
R_1 m_0 - R_2 \leq  c(a, f) \leq R_1 m_0 + R_2,
\]
and hence $c (a, f)$ is finite.
The rest of the proof in
\cite{Oh-act}
of the finiteness and continuity of spectral invariants
can be carried over to our case in a direct fashion.
In particular, the continuity follows from the following
$C^0{\hbox{\rm -estimate}}$ \cite{Oh-act}, \cite{Sch-1}.
For any Hamiltonians $F^\prime,
F^{{\prime\prime}}\in \cF$ and $a \in QH_{ev}(M)$ one has
\begin{eqnarray}
\label{eqn-c-0-estimate}
\int_{S^1} -\max_M\, (F_t^\prime  - F_t^{{\prime\prime}}) \, dt \leq
& c (a, \tilde{\psi}_{F^\prime}) - c (a, \tilde{\psi}_{F^{{\prime\prime}}})
\leq &\nonumber\\
\leq &
\displaystyle{\int_{S^1} -\min_M\, (F_t^\prime  - F_t^{{\prime\prime}}) \,
dt.} &
\end{eqnarray}

The numbers $c(a,f)$ are
called {\it spectral invariants} of an element $f \in \tG$.
One can show that they all lie in
the action spectrum $\spec(f)$ \cite{Oh-act},\cite{Oh-new},
\cite{Sch-1}
(this is obviously true for a generic Hamiltonian from the original
definition (\ref{eqn-D-def}) of
$d$) and
persist under conjugations in $\tG$
(see \cite{Oh-act}, \cite{Oh-new},
\cite{Sch-1}, cf.
\cite{Entov}).
Spectral invariants will play a crucial role below
in the construction of the
Calabi quasimorphism.

\subsubsection{Spectral invariants as
characteristic exponents} \label{charexp}

A function $\chi : V \to \R \cup {-\infty}$
on a vector space $V$ over $\C$
is called
{\it a characteristic exponent}
if
\begin{itemize}
\item $\chi (v) \in \R $ for all
non-zero $v \in V$ and $\chi(0) = -\infty$;
\item $\chi (\delta \cdot v) = \chi (v) $ for
every non-zero $\delta \in \C$
and $v \in V$;
\item $\chi (v_1 + v_2) \leq
\max(\chi (v_1),\chi_(v_2))$ for all $v_1,v_2 \in V$.
\end{itemize}
This notion (which appears in the theory
of Lyapunov exponents in Dynamical
Systems, see e.g. \cite{DS})
is relevant in our study of spectral invariants.
It is a straightforward consequence of the definitions
that for a given $f \in \tG$ the function
$$c(\cdot,f): QH_{ev}(M) \to \R,\; a \mapsto c(a,f)$$
is a characteristic exponent on $QH_{ev}(M)$.
Starting from this observation, one can apply various known facts
about characteristic exponents to the spectral invariants.
For instance, for every $f \in \tG$ and $m \in \Z$ the set
$$\{c(a,f)\;|\; a \in QH_{2m}(M)\}$$
has at most $\dim_{\C} QH_{2m}(M)$ distinct elements.
The well known
fact which will be used below is as follows. We formulate it
in the language of characteristic exponents.

\begin{prop}
\label{char}
Let $\chi: V \to \R$ be a characteristic exponent.
Assume that $\chi(v_1) < \chi (v_2)$.
Then $\chi (v_1 + v_2) = \chi (v_2)$.
\end{prop}

\bigskip
\noindent
{\bf Proof:}
By definition, $\chi (v_1 + v_2) \leq \chi (v_2)$.
Assume on the contrary that $\chi (v_1 + v_2) < \chi (v_2)$.
Then (using that $\chi(-v_1) = \chi(v_1)$) we have
$$\chi (v_2) = \chi (-v_1 + (v_1 +
v_2)) \leq \max(\chi (v_1),\chi(v_1+v_2)) < \chi(v_2),$$
a contradiction. Hence $\chi(v_1+v_2) = \chi(v_2)$.
$\square$

\subsubsection{Spectral invariants of the identity}\prlabel{2.6.D}

Define a function $\nu : QH_{ev}(M) \to \Z$ as follows.
For a non-zero element $a = \sum A_j s^j \in QH_{ev}(M)$
set $\nu (a)$ to be the maximal $j$ such that $A_j \neq 0$.
We claim that
\begin{equation}
\label{eqn-c-id}
c(a, \id_\tG) = \Omega \nu (a),
\end{equation}
where $\id_\tG$ stands for the identity in $\tG$.

Indeed, according to
\cite{Oh-act} and
\cite{Sch-1},
$c (A, \id_\tG ) = 0$
for any singular homology class
$A \neq 0$. Suppose now that $a \neq 0$ is an arbitrary
quantum homology class. In view of formula (\ref{eqn-c-max}) above
it suffices to prove the claim assuming that $a$ is homogeneous
in the sense of the grading.
In this case  $a$ is given by
a {\it finite} sum of the form
$a = \sum_m A_{m} s^m$, where
$A_{m}\in H_{ev} (M)$.

Using the
$k{\hbox{\rm -linearity}}$
of the identification  between
$V_{\infty} (f)$ and $QH_{ev}(M)$
and
formula (\ref{act2}) one easily gets that
$$c(sb,f) = c(b,f) +
\Omega \;\text{for}\;\text{all}\; b \in QH_{ev}(M),\; f \in \tG.$$
Therefore
$c (A_{m} s^m, \id_\tG) = \Omega m + c (A_{m}, \id_\tG) = \Omega m$
for any
$m$
such that
$A_{m} \neq 0$.
In view of Proposition \ref{char}
$$c(a,\id_\tG) = \max_{m:A_{m}\neq 0} \Omega m = \Omega \nu (a),$$
and the claim follows (cf. \cite{Oh-new}).

In fact, given an arbitrary $a\in QH_{ev} (M)$
one can calculate $c(a, f)$ not only for
$f=\id_\tG$ but for any $f\in \pi_1 (G)\subset \tG$, i.e. for any
lift of $\id\in G$ to $\tG$ -- see Section~\ref{sect-c-on-loops}.

\subsubsection{Pair-of-pants product} \prlabel{2.6.E}
There exists
a so-called {\it pair-of-pants product} \cite{PSS}
$$V_{\alpha}(f) \times V_{\beta} (g) \to V_{\alpha+\beta}(fg)
,\;\; (v,w) \mapsto v*_{PP} w.$$ It agrees with the quantum
product, namely $$i_{\alpha+\beta}(v*_{PP} w) = i_{\alpha}(v) *
i_{\beta}(w).$$ This immediately yields the following triangle
inequality for spectral invariants:
$$c(a*b,fg) \leq c(a,f) + c(b,g),\;\; a,b \in QH_{ev}(M).$$

\subsubsection{Natural projections}\prlabel{2.6.F}
The natural projection of Floer complexes $C_{\infty}(F) \to
C_{\infty}(F)/C_{\alpha}(F)$
 induces a homomorphism
\[
\pi_{\alpha} : QH_{ev}(M) = V_{\infty}(f) \to V^{\alpha} (f).
\]
between their homology. The homological exact sequence yields
$\text{Kernel}\; \pi_{\alpha} = \text{Image}\; i_{\alpha}.$ In the
case we wish to emphasize the dependence of $\pi_{\alpha}$ on the
element $f$ we will write $\pi_{\alpha}\{f\}$.

\subsubsection{Poincar{\'e} duality}\prlabel{2.6.G}

Each critical point
$y$
of
$\cA_F$
of Conley-Zehnder index
${\textit ind}\, (y)$
is also a critical point of
$\cA_{-F}$
of Conley-Zehnder index $2n -{\textit ind}\, (y)$.
Moreover there exists a non-degenerate intersection pairing
between the Floer complex associated to
$F$
and the Floer complex of
$-F$
leading to the Poincar{\'e} duality in the Floer homology theory
similarly to the situation in the classical Morse homology theory.
One can show that for every $\alpha \in \R$
the space $V_{\alpha}(f^{-1})$ is canonically
isomorphic to $\text{Hom}~(V^{-\alpha}(f),\C)$.
This isomorphism gives rise to a non-degenerate pairing
$L: V_{\alpha}(f^{-1}) \times V^{-\alpha} (f) \to \C$
which agrees with the intersection pairing $\Pi$
on the quantum homology
(see equation (\ref{3.2.A}) of Section \ref{2.3}):
$$\Pi(i_{\alpha}\{f^{-1}\} a, b) = L(a, \pi_{-\alpha}\{f\}b)$$
for every $a \in V_{\alpha}(f^{-1})$ and $b\in QH_{ev}(M)$.
These statements can be extracted from \cite{PSS}.

\subsection{Comparing spectral invariants of $f$ and $f^{-1}$}
\prlabel{2.7}

For an element $b \in QH_{ev}(M)\setminus\{0\}$
denote by $\Upsilon(b)$ the set of
all $a \in QH_{ev}(M)$ with $\Pi(a,b) \neq 0$.

\begin{lemma}
\prlabel{2.7.A}
$$c(b,f) = -\inf_{a \in \Upsilon(b)} c(a,f^{-1})$$
for all $b \in QH_{ev}\setminus\{0\}$ and $f \in \tG$.
\end{lemma}

\bigskip
\noindent
{\bf Proof:}
Set $$\delta =
\inf_{a \in \Upsilon(b)} c(a,f^{-1}).$$
The proof is divided into two steps.

\medskip
\noindent
1) Take arbitrary $\epsilon > 0$ and set $\alpha = \epsilon - c(b,f)$.
Then $b \notin \text{Image}\;i_{-\alpha}\{f\}$, so by \ref{2.6.F}
$$w:=\pi_{-\alpha}\{f\}b \neq 0.$$
Since the pairing $L$ is non-degenerate, there exists $v \in
V_{\alpha}(f^{-1})$ such that \break $L(v,w) \neq 0$. Put $a_0 =
i_{\alpha}\{f^{-1}\}v$. By \ref{2.6.G} we have $L(v,w) =
\Pi(a_0,b) \neq 0$. We see that $c(a_0,f^{-1}) \leq \alpha$, so
$\delta \leq \alpha = \epsilon - c(b,f)$. Since this inequality
holds for every $\epsilon > 0$ we conclude that $\delta \leq
-c(b,f)$.

\medskip
\noindent
2) Take arbitrary $\epsilon > 0$ and set $\alpha = -c(b,f) - \epsilon$.
Then $b \in \text{Image}\; i_{-\alpha}\{f\}$, so by \ref{2.6.F}
$\pi_{-\alpha}\{f\}b = 0$. Assume that there exists $a \in \Upsilon(b)$
such that $c(a,f^{-1})< \alpha$. Then
$a \in \text{Image}\;i_{\beta}\{f^{-1}\}$ for some $\beta < \alpha$,
and hence
$a \in \text{Image}\;i_{\alpha}\{f^{-1}\}$ in view of \ref{2.6.B}.
Take $v \in V_{\alpha}(f^{-1})$ so that $a = i_{\alpha}\{f^{-1}\}v$.
By \ref{2.6.G}
$$\Pi(a,b) = L(v,\pi_{-\alpha}\{f\}b) = 0,$$
and we get a contradiction with the assumption $\Pi(a,b)\neq 0$.
Hence
\break
$c(a,f^{-1}) \geq \alpha$ for every $a \in \Upsilon(b)$, so
$\delta \geq \alpha= -c(b,f) -\epsilon$. Since this is true
for every $\epsilon > 0$ we get that $\delta \geq -c(b,f)$.
Combining this with the inequality proved in Step 1 we conclude
that $\delta = -c(b,f)$ as required.

\hfill$\square$

\section{Constructing the quasimorphism} \prlabel{3}

\subsection{Spectral numbers define a quasimorphism $r$ on $\tG$} \prlabel{3.1}
Suppose that the algebra $Q = QH_{ev}(M)$ is semi-simple, and let
$Q= Q_1 \oplus ... \oplus Q_d$ be its decomposition into the
direct sum of fields. Denote by $e$ the unity of $Q_1$.

\begin{thm}\prlabel{3.1.A}
The function
$$r: \tG \to \R,\;\; f \mapsto c(e,f)$$
is a quasimorphism.
\end{thm}

\bigskip
The proof is given in Section \ref{3.4} below. For the proof
we need the following lemma.

\subsection{A lemma from non-Archimedian geometry}
\prlabel{3.3}

Let $\nu: Q \to \Z$ be the function introduced in \ref{2.6.D}.

\begin{lemma}
\prlabel{3.3.A}
There exists $R>0$ such that $\nu(b) + \nu(b^{-1}) \leq R$
for every $b \in Q_1 \setminus \{0\}$.
\end{lemma}

\bigskip
We are grateful to V.~Berkovich for explaining to us the proof.
The reader is referred to \cite{Gouvea} for preliminaries
on non-Archimedian geometry.

\bigskip
\noindent
{\bf Proof:} For $\varsigma \in k$ set $|\varsigma|
= \exp \nu(\varsigma)$.
Then $|\;|$ is a non-Archimedian absolute
value on $k$, and the field $k$
is complete with
respect to $|\;|$. For $b \in Q_1$ put $||b|| = \exp \nu (b)$.
Then $||\;||$ is a
norm on $Q_1$, where $Q_1$ is considered as a vector space
over $k$. Since the field $Q_1$ is a finite extension of $k$,
the absolute value $|\;|$ extends
to an absolute value $|||\;|||$ on $Q_1$.
Furthermore,
all norms on a finite-dimensional space over $k$ are equivalent.
Thus there exists $\delta > 0$ so that
$$||b|| \leq \delta \cdot |||b|||$$
for every $b \in Q_1$. Therefore for $b\neq 0$
$$||b||\cdot ||b^{-1}|| \leq
\delta^2\cdot |||b|||\cdot|||b^{-1}|||= \delta^2,$$
where the last
equality follows from the definition of the absolute value.
Therefore $\nu(b)+\nu(b^{-1}) \leq R$ with $R = 2\log \delta$.
This completes the proof.
\hfill $\square$

\subsection{Proof of Theorem~\ref{3.1.A}}
\prlabel{3.4}

Note that $e*e=e$. By the inequality in Section~\ref{2.6.E} we have
\begin{equation}
\prlabel{3.4.A}
c(e,fg)= c(e*e,fg) \leq c(e,f)+c(e,g).
\end{equation}
Similarly,
$$c(e,fg) \geq c(e,f) - c(e,g^{-1}).$$
Applying Lemma~\ref{2.7.A}
we get that
\begin{equation}
\prlabel{3.4.B}
c(e,fg) \geq c(e,f) + \inf_{b:\Pi(b,e)\neq 0} c(b,g).
\end{equation}
Our next goal is to find a lower bound for $c(b,g)$ provided
$\Pi(b,e) \neq 0$. Write $b=b_1+...+b_d$ where $b_i \in Q_i$ for
all $i=1,...,d$. Then formula (\ref{3.2.A}) of Section~\ref{2.3}
yields $$\Pi(b,e) = \tau\Delta(b*e,[M]) = \tau\Delta(b_1,[M]) \neq
0.$$ This immediately implies that $b_1 \neq 0$ and $\nu (b_1)
\geq 0$. Therefore $b_1$ is invertible in $Q_1$ and $\nu
(b_1^{-1}) \leq R$, where $R$ is the constant from
Lemma~\ref{3.3.A}. Applying \ref{2.6.E} we obtain $$c(b,g) \geq
c(b*e,g) - c(e,\id_\tG) = c(b_1*e,g)- c(e,\id_\tG)$$ $$ \geq
c(e,g) -c(b_1^{-1},\id_\tG) -c(e,\id_\tG).$$ Using \ref{2.6.D} we
get that $$c(b_1^{-1},\id_\tG) = \Omega\nu(b_1^{-1}) \leq \Omega
R,$$ and therefore $$c(b,g) \geq c(e,g) - \Omega R -
c(e,\id_\tG).$$ Substituting this into (\ref{3.4.B}) we see that
$$c(e,fg) \geq c(e,f)+c(e,g) - \hbox{\it const}.$$ Together with
(\ref{3.4.A}) this proves that the map  $r$ which sends $f$ to
$c(e,f)$ is a quasimorphism. \hfill$\square$

\subsection{Building a Calabi quasimor\-phism $\tmu$ from $r$}
\prlabel{3.5}

Consider now a homogeneous quasimorphism $\tmu:\tG \to \R$
given by
\begin{equation}
\label{eqn-def-mu} \tmu (\tf) =-\hbox{\rm vol}\, (M)\,\cdot
\lim_{m \to \infty}\frac{r(\tf^m)}{m},
\end{equation}
where $r(\tf)= c(e,\tf) $ and $$\hbox{\rm vol}\, (M)\,=\int_M
\omega^n.$$ Let $\cD$ be the class of all displaceable open
subsets of $M$ as in (\ref{1.1.C}).

\begin{prop}
\prlabel{3.5.A}
The restriction of $\tmu$ on $\tG _U$ coincides
with the Calabi homomorphism $\tCa _U$ for every
$U \in \cD$.
\end{prop}

\bigskip
\noindent {\bf Proof:} We follow closely the work by Yaron
Ostrover \cite{Ostr}. Take an open subset $U \in \cD$. By
definition there exists a Hamiltonian diffeomorphism $h \in G$
which displaces $U$: $$h(U) \cap \text{Closure}\;(U) =
\emptyset.$$ Fix any lift $\tih$ of $h$ to $\tG$. Let $F: M \times
\R \to \R$ be a Hamiltonian function which is 1-periodic in time
and satisfies $F(x,t) = 0$ for all $t \in \R, x \in M \setminus
U$. Write $f_t$ for the corresponding Hamiltonian flow, and $\tf
_t$ for its lift to $\tG$. Put $\tf = \tf_1$ and note that the
periodicity of $F$ in $t$ yields $\tf_m = \tf^m$ for all $m \in
\Z$. Put $$F'(x,t) = F(x,t) -{\{\hbox{\rm vol}\, (M)\} }^{-1}\,
\cdot \int_M F(x,t)\, \omega^n.$$ Note that $F'(x,t)$ generates
the same flow $f_t$ and satisfies the normalization condition
$\int_M F'(x,t)\, \omega^n = 0$ which enters the definition of the
action functional (see Section~\ref{2.4}). Consider the family
$\tih\tf_t,\; t \in \R$. Since $hU \cap U = \emptyset$ and $f_t(U)
= U$ the fixed point set of $hf_t$ coincides with the fixed point
set of $h$ for every $t$, and hence lies outside $U$. Moreover,
for every $t$ and every $x \in M\setminus U$ we have $F'(x,t)
\equiv u (t)$ where $$u(t) = -{\{\hbox{\rm vol}\, (M)\} }^{-1}\,
\cdot \int_M F(x,t)\, \omega^n.$$ Using this one can calculate the
action spectrum of $hf_t$: $$\spec (\tih \tf _t) = \spec (\tih) +
w (t),\;\; \text{where}\; w(t) := \int_0^t u (z) dz.$$ Consider
the function $\psi (t) = r(\tih\tf_t) = c(e,\tih\tf_t)$. It is
continuous and takes values in $\spec (\tih\tf_t)$ (see
\ref{2.6.C}). Since $\spec (\tih)$ is a closed nowhere dense
subset of $\R$ we conclude that there exists $s_0 \in \spec (h_0)$
such that $\psi (t) = s_0 + w(t)$. Hence $r(\tih\tf^m) = s_0 +
w(m)$. Using that $r$ is a quasimorphism we calculate $$\tmu (\tf)
= -\hbox{\rm vol}\, (M)\,\cdot \lim_{m\to +\infty} \frac{r(\tih
\tf^m)}{m} = -\hbox{\rm vol}\, (M)\, \cdot \lim_{m\to +\infty}
\frac{w (m)}{m}$$ $$= \int_0^1 dt \int_M F(x,t)\, \omega^n =
\tCa_U (\tf).$$ This completes the proof. \hfill$\square$

\subsection{From $\tG$ to $G$}\prlabel{3.6}

\begin{prop}
\prlabel{3.6.A}

\noindent
Suppose that $\pi_1(G)$ is finite. Then the quasimorphism
$\tmu: \tG \to \R$ descends to a homogeneous Calabi quasimorphism
$\mu : G \to \R$.

\end{prop}

\bigskip
\noindent
{\bf Proof:} The fundamental group $\pi_1 (G)$ is the
kernel of the natural projection $\tG \to G$. Note that
$\pi_1 (G)$ lies in the center of $\tG$. Then
$$\tmu (\phi\tf) = \tmu(\phi) + \tmu(\tf)$$
for every $\phi \in \pi_1(G)$ and $\tf \in \tG$.
This is true since for any
$m$
the quantities
$m\tmu (\phi\tf) = \tmu ((\phi\tf)^m)$
and
$\tmu (\phi^m) +\tmu (\tf^m) = m (\tmu(\phi) + \tmu(\tf) )$,
which are homogeneous with respect to $m$,
differ by a constant which is independent of
$m$.
Since $\pi_1(G)$ is finite $\tmu$ vanishes on $\pi_1(G)$.
Hence
$\tmu (\phi\tf) = \tmu(\tf)$, so $\tmu$ descends to a function
$\mu: G \to \R$. Using  Proposition \ref{3.5.A} one readily checks that
$\mu$ is the required quasimorphism.
\hfill $\square$

\subsection{Hofer metric and
the continuity properties of $\tmu$ and $\mu$}
\prlabel{3.7}

As before we write $\cF$ for the space of all normalized
time-dependent Hamiltonians on $M$. This space is equipped
with a $C^0$-norm given by formula (\ref{conorm}).
We denote by ${\tilde \psi}_F$ the element of $\tG$ generated
by a Hamiltonian $F \in \cF$.

The Hofer metric on $G$ (see Section~\ref{1.3.BB}) can be lifted
to a bi-invariant pseudo-metric $\tilde{\rho}$ on $\tG$
defined as follows. Let $\psi\in G$, let
$\tilde{\psi}\in
\tG$
be its lift, let $\id$ be the identity in $G$ and $\id_{\tG}$ the
identity in $\tG$. Then
\[
\tilde{\rho} (\id_\tG, \tilde{\psi}) = \inf_F \int_{S^1}
\|F_t\|_{C^0} \, dt,
\]
where the infimum is taken over all $F \in \cF$ such that $\tilde{\psi}
=\tilde{\psi}_F$.
In particular,
\begin{equation}
\label{eqn-rho-rho-tilde}
\rho (\id, \psi) = \inf_F \tilde{\rho} (\id_\tG, \tilde{\psi}_F),
\end{equation}
where the infimum is taken over all Hamiltonians $F$ generating
$\psi$ (or, equivalently, over all lifts $\tilde{\psi}_F$ of $\psi$
to $\tG$).

We will prove now that $\tmu$ is a continuous function on $\tG$ and
is Lipschitz with respect to $\tilde{\rho}$.

\begin{prop}
\label{prop-tmu-lipschitz}
Let $\tmu: \tG\to\R$ be the Calabi quasimorphism constructed
above.
Then for any Hamiltonians
${H^\prime}, {H^{\prime\prime}}$
\begin{eqnarray}
\label{eqn-tmu-lipschitz} | \tmu (\tilde{\psi}_{H^\prime}) -\tmu
(\tilde{\psi}_{H^{\prime\prime}})| \leq & | \hbox{\rm vol}\,
(M)|\, \cdot \tilde{\rho}
(\tilde{\psi}_{H^\prime},\tilde{\psi}_{H^{\prime\prime}})\leq
&\nonumber\\ \leq & | \hbox{\rm vol}\, (M)|\, \cdot \int_{S^1}
\;\| {H_t^\prime} - {H_t^{\prime\prime}}\|_{C^0}\; dt, &
\end{eqnarray}
and therefore $\tmu: \tG\to\R$
is continuous.
\end{prop}

\bigskip
Proposition~\ref{prop-tmu-lipschitz}
together with (\ref{eqn-rho-rho-tilde})
immediately lead to the following corollary.

\begin{cor}
\label{cor-mu-lipschitz}
If the quasimorphism $\tmu: \tG\to\R$ constructed above
descends to a quasimorphism  $\mu: G\to\R$ then
for any Hamiltonians
${H^\prime}, {H^{\prime\prime}}$
\begin{eqnarray}
\label{eqn-mu-lipschitz} | \mu (\psi_{H^\prime}) -\mu
({\psi}_{H^{\prime\prime}}) |\leq & | \hbox{\rm vol}\, (M)| \,
\cdot {\rho} ({\psi}_{H^\prime},{\psi}_{H^{\prime\prime}})\leq
&\nonumber\\ \leq & |\hbox{\rm vol}\, (M)|  \, \cdot \int_{S^1}
\;\| {H_t^\prime} - {H_t^{\prime\prime}}\|_{C^0} \; dt,&
\end{eqnarray}
and therefore $\mu: G\to\R$ is a continuous function.
\end{cor}

\bigskip
\noindent
{\bf Proof of Proposition~\ref{prop-tmu-lipschitz}:}
Suppose that ${H^\prime}$ and ${H^{\prime\prime}}$ generate Hamiltonian flows
$\{ f_t\}$ and $\{ g_t\}$ which give rise to the elements
$f:=\tilde{\psi}_{H^\prime}$ and $g: =\tilde{\psi}_{H^{\prime\prime}}$.
Recall that the Hamiltonian ${H^\prime}\sharp {H^{\prime\prime}} (x, t) :=
{H^\prime} (x, t) + {H^{\prime\prime}} (g^{-1}_t (x), t)$ generates the flow
$\{ f_t g_t\}$
and the Hamiltonian $\overline{{H^\prime}} (x,t) := - {H^\prime} (f_t (x), t)$
generates
the flow $\{ f_t^{-1} \}$. Set $H : = \overline{{H^\prime}}\sharp {H^{\prime\prime}}$.
Thus
\[
\tilde{\psi}_{H^{\prime}}^{-1}\tilde{\psi}_{H^{\prime\prime}} =
\tilde{\psi}_H,
\]
where $H (x,t) = - {H^\prime} (f_t (x), t) + {H^{\prime\prime}}
(f_t (x), t)$. Observe that for each $t$
\begin{equation}
\label{eqn-C-0-norm} \| H_t\|_{C^0} = \| {H_t^\prime} -
{H_t^{\prime\prime}}\|_{C^0}.
\end{equation}
Since $\tilde{\rho}$ is bi-invariant we have
\[
\tilde{\rho} (f, g) = \tilde{\rho} (\id_\tG, f^{-1} g) \leq
\int_{S^1} \| H_t\|_{C^0} \, dt \; .
\]
Combining this with (\ref{eqn-C-0-norm})
we get the second inequality in
(\ref{eqn-tmu-lipschitz}).

Now let us prove the first inequality in
(\ref{eqn-tmu-lipschitz}).
Indeed, according to
(\ref{eqn-c-0-estimate}), for any $a\in QH_\ast (M)$
\[
| c (a, f) - c (a, g) |  \leq \int_{S^1} \| H^\prime_t -
H^{\prime\prime}_t \|_{C^0} \, dt = \int_{S^1} \|H_t\|_{C^0}\; dt.
\]

This inequality is true for {\it any} ${H^\prime}$ and ${H^{\prime\prime}}$ generating, respectively,
$f$ and $g$, while its left-hand side
depends only on the elements $f, g\in\tG$
and not on the Hamiltonians that generate them. Thus taking in the right-hand side
the infimum over all ${H^\prime}$ and ${H^{\prime\prime}}$ generating, respectively,
$f$ and $g$, we obtain
\[
| c (a, f) - c (a, g) |  \leq
\tilde{\rho} (\id_\tG, f^{-1} g )
\]
and hence
\begin{equation}
\label{eqn-c-f-m-g-m-tilde-rho}
| c (a, f^m) - c (a, g^m) |  \leq
\tilde{\rho} (\id_\tG, f^{-m} g^m )
\end{equation}

Now observe that
\[
f^{-m} g^m = \prod_{i=0}^{m-1} g^{-i} (f^{-1} g) g^i.
\]
Thus
\begin{equation}
\label{eqn-group-identity}
\tilde{\rho} (\id_\tG, f^{-m} g^m )\leq
\sum_{i=0}^{m-1}
\tilde{\rho} (\id_\tG, g^{-i} (f^{-1} g) g^i )
\leq m \tilde{\rho} (\id_\tG, f^{-1} g ),
\end{equation}
where the last inequality holds because
$\tilde{\rho}$ is bi-invariant.
Combining (\ref{eqn-c-f-m-g-m-tilde-rho}) with
(\ref{eqn-group-identity}) we see that
\begin{equation}
\label{eqn-c-tilde-rho-fixed-m}
\frac{1}{m} | c (a, f^m) - c (a, g^m) |  \leq
\tilde{\rho} (\id_\tG, f^{-1} g ) = \tilde{\rho} (f, g).
\end{equation}
Now take $a$ to be the unit element $e$ of the field $Q_1$
involved in the definition of the quasimorphism $r = c (e, \cdot):
\tG\to\R$ (see Section~\ref{3.1}) and recall that, according to
its definition, $\tmu (f) = - \hbox{\rm vol}\, (M) \lim_{m\to
+\infty}  c (e, f^m)/m$ (see Section~\ref{3.5}). Together with
(\ref{eqn-c-tilde-rho-fixed-m}) this yields
\[
|\tmu (f) - \tmu (g) | \leq | \hbox{\rm vol}\, (M) |\, \cdot
\tilde{\rho} (f, g),
\]
proving the first inequality in (\ref{eqn-tmu-lipschitz}). The
proposition is proven.
\hfill$\square$

\subsection{Proofs of Theorems~\ref{1.1.D} (for $M =
S^2, S^2\times S^2, \C P^2$)
and \ref{1.1.E}}
\prlabel{3.8}

\medskip
\noindent
{\bf Proof of Theorem~\ref{1.1.E}:} According to
Theorem~\ref{3.1.A} and Proposition \ref{3.5.A}, the function
$\tmu : \tG\to\R$ constructed above is a homogeneous Calabi quasimorphism.
In view of Proposition~\ref{prop-tmu-lipschitz} it is continuous
(and even Lipschitz with respect to the Hofer pseudo-metric on
$\tG$).
This proves Theorem~\ref{1.1.E}.
\hfill$\square$

\bigskip
\noindent
{\bf Proof of Theorem~\ref{1.1.D} for $M = S^2, S^2\times S^2, \C P^2$:}
Let $(M,\omega)$ be one of the manifolds $S^2, S^2\times S^2, \C P^2$.
Then $(M,\omega)$ is
spherically monotone.
The quantum homology algebra $QH_{ev}(M)$
is semi-simple (see
Section~\ref{sec-qh-examples}). Thus Theorem~\ref{1.1.E}
gives us a continuous homogeneous Calabi quasimorphism
$\tmu:\tG\to\R$. The fundamental group $\pi_1(G)$ is finite for
$M = S^2, S^2\times S^2, \C P^2$ (see \cite{Gro-pshc}).
Therefore, according to
Proposition~\ref{3.6.A}, $\tmu$ descends to a homogeneous Calabi quasimorphism
on $\mu: G\to\R$. In view of Corollary~\ref{cor-mu-lipschitz},
$\mu$ is continuous
(and Lipschitz with respect to the Hofer metric on
$G$). The theorem is proven.
\hfill$\square$

\section{Spectral invariants and Hamiltonian loops}
\label{sect-c-on-loops}

A homogeneous quasimorphism on an abelian group is always a
homomorphism (the proof of this simple fact is actually contained
in the proof of Proposition~\ref{3.6.A} above). Thus the
restriction of the quasimorphism $\tmu: \tG\to\R$ constructed
above on the abelian subgroup $\pi_1 (G)\subset \tG$ is a
homomorphism. In this section we obtain a formula for
this homomorphism in terms of the Seidel action of $\pi_1(G)$ on
the quantum homology of $M$ \cite{Se,Lal-McD-Pol}. As an
application we show that this homomorphism vanishes when $M$ is
the projective space $\C P^n$  endowed with the Fubini-Study form
and thus complete the proof of Theorem \ref{1.1.D}. The results of
this section were communicated to us by Paul Seidel.

\subsection{Preliminaries on the Seidel action}

\subsubsection{Hamiltonian fibrations over $S^2$}

There exists a one-to-one correspondence between homotopy classes
of loops in $G$ and isomorphism classes of Hamiltonian fibrations
over the 2-sphere $S^2$ with the fiber $(M^{2n},\omega)$, see
\cite{McD-Sal-sympl-top, Pol-book}. We denote by $\pi:E_\gamma \to
S^2$ the fibration associated to a loop $\gamma$. An important
invariant of such a fibration is its {\it coupling class}  $W\in
H^2 (E_\gamma, \R)$ which is uniquely defined by the following
conditions: the restriction of $W$ to each fiber coincides with
the class of the symplectic form, and its top power $W^{n+1}$
vanishes.

Denote by $T^{vert} E_\gamma$
the vector bundle over $E_\gamma$ formed by all tangent spaces of
the fibers of $\pi$ and by $c_1^{vert}$ the first Chern class of this
bundle.

Take a positively oriented complex structure $j$ on $S^2$ and an
almost complex structure $\hJ$ on $E_\gamma$ whose restriction on
each fiber is compatible with the symplectic form on it and such that
the projection $\pi$ is a $(\hat{J}, j){\hbox{\rm -holomorphic}}$ map
(see \cite{Se}).

Two $(j,\hJ){\hbox{\rm -holomorphic}}$
sections $v_1, v_2$ of $\pi: E_\gamma\to S^2$  are said
to be equivalent if $W ([v_1 (S^2)] ) = W ([v_2 (S^2)] )$.
Since $M$ is assumed to be spherically monotone this condition is
equivalent to
$c_1^{vert} ([v_1 (S^2)]) = c_1^{vert} ([v_2 (S^2)] )$.
Denote by
${\cS}_\gamma$
the set of all such equivalence classes --
it is an
affine space modeled on $\bar{\pi}_2 (M) \cong \Z$.
According to the definition, the maps
$v \mapsto W ([v (S^2)] )$ and
$v \mapsto c_1^{vert} ([v (S^2)] )$ give rise to some
correctly defined functions on ${\cS}_\gamma$. From this moment
on,
given an equivalence class
$\sigma\in {\cS}_\gamma$, we will denote by $W (\sigma)$, $c_1^{vert} (\sigma)$
the values of those functions
on $\sigma$.

Recall that in our notation $S$ is the
positive generator of $\bar{\pi}_2
(M)$
and
$\Omega : = (\omega, S) > 0$, $N := (c_1 (M), S)$.
Thus for a class
$\sigma\in {\cS}_\gamma$ a sum $\sigma + mS$, $m\in\Z$,
stands for another class in
${\cS}_\gamma$ so that:
\begin{equation}
\label{eqn-u-c-1-sigma-S}
W (\sigma + mS) = W (\sigma) +m\Omega,\
c_1^{vert} (\sigma + mS) = c_1^{vert}
(\sigma) + mN.
\end{equation}

\subsubsection{Gromov-Witten invariants revisited}

Given $\sigma\in {\cS}_\gamma$ and homology classes $A, B, C\in
H_{ev} (M)$ let $GW_\sigma (A, B, C)$ denote the Gromov-Witten
number (cf. Section \ref{2.2}) , defined as follows. In the fibers
of $\pi$ over $0,1,\infty\in S^2$ pick the cycles $\widehat{A},
\widehat{B}, \widehat{C}$ realizing, respectively, the homology
classes $A, B, C$. Consider the sections from the class $\sigma$
whose intersection with the fibers over $0,1,\infty$ belongs,
respectively, to $\widehat{A}, \widehat{B}, \widehat{C}$. If there
is a finite number of such sections count them with appropriate
signs and set $GW_\sigma (A, B, C)$ equal to the result. Otherwise
set $GW_\sigma (A, B, C)$ to be zero. The resulting number does
not depend on the choice of cycles
$\widehat{A},\widehat{B},\widehat{C}$ -- for details see
\cite{Se}.

\subsubsection{Hamiltonian loops}

Consider the space $\La$ of all smooth contractible loops in $M$
(i.e. smooth maps from $S^1$ to $M$). Let $\tLa$ be the cover of
$\La$ introduced in Section \ref{2.4}. Its elements are
equivalence classes of pairs $(x,u)$, where $x \in \La$,  $u$ is
an oriented disk spanning $x$ in $M$, and the equivalence relation
is defined as follows: $(x_1,u_1) \sim (x_2,u_2)$ iff $x_1 = x_2$
and the 2-sphere $u_1 \# (-u_2)$ vanishes in $\bar{\pi}_2 (M)$.

The group $\cG$ of all (smooth) identity-based loops in $G$ acts on
$\La$: if $\gamma = \{ g_t\}\in\cG$ then the action $T_{\gamma} :
\La\to\La$ is defined as
\[
T_{\gamma} \{ x_t\} = \{ g_t (x_t)\}.
\]
This map can be lifted (not uniquely!) to a map on $\tLa$. In fact
there is a one-to-one correspondence between lifts of $T_{\gamma}$
and classes of sections $\sigma \in \cS_{\gamma}$. We denote the
lift corresponding to $\sigma$ by ${\tilde T}_{\gamma,\sigma}$.

Suppose that  Hamiltonian loop $\gamma=\{ g_t\}\in \cG$ is
generated by a normalized Hamiltonian $K: M\times S^1\to\R$, $K\in
\cF$, and let $H\in \cF$. Consider the action functional $\cA_H$
on $\tLa$ (see Section~\ref{2.4}). The following formula
 (see
\cite{Se, Lal-McD-Pol}) is crucial for our purposes:
\begin{equation}\label{vsp1}
({\tilde T}_{\gamma,\sigma}^*)^{-1} \cA_H -\cA_{K\sharp H} =
-W(\sigma).
\end{equation}
In particular, the function in the left hand side is constant on
$\tLa$.

The action of ${\tilde T}_{\gamma,\sigma}$
on $\tLa$ defines an isomorphism,
which we will denote by
$\iota$,
between the Floer homology of $H$ and the Floer homology of
$K\sharp H$ \cite{Se}.
According to Seidel's theorem \cite{Se}, under the identification
of Floer and quantum homology (see Section~\ref{2.6.A}) this
isomorphism between Floer homology groups
corresponds to the
multiplication by some class $\Psi$ in the quantum homology of $M$.
This class is defined as follows:
\begin{equation}
\label{vsp3}
\Psi =
\sum_{m\in\Z} A_{\sigma + mS} s^{-m},
\end{equation}
where $A_{\sigma +mS} \in H_{ev} (M)$ is uniquely determined
 by the condition
\[
A_{\sigma +mS}\circ_M C = GW_{\sigma + mS} ([M], [M], C)
\]
for any
$C\in H_{ev} (M)$.

Note also that in view of
(\ref{vsp1}), the isomorphism $\iota$
shifts the filtration of the Floer homology groups by $W(\sigma)$:
\begin{equation}\label{vsp2}
\iota : V_\alpha (H)\to V_{\alpha + W(\sigma)} (K\sharp H)
\end{equation}
for any $\alpha\in\R$.

\subsubsection{Extending the field} \label{subsec-field}

In what follows it would be convenient to work with
the extension $\bar{k}$ of the field $k$, where $\bar{k}$
is formed by semi-infinite sums $\sum_{\alpha\in\R} z_\alpha
s^\alpha$, $z_\alpha\in\C$, satisfying the condition that for any
$\alpha_0\in\R$ there is only a finite number of terms with
$z_\alpha\neq 0$, $\alpha\geq \alpha_0$, in the sum. Define
$${\oQ} := H_{ev} (M)\otimes_\C \bar{k} =
QH_{ev} (M)\otimes_k \bar{k}.$$
Naturally, $QH_{ev} (M)\subset
\oQ$. The function $\nu$ on $QH_{ev}(M)$ defined in
Section \ref{2.6.D}
extends to a function on $\oQ$
which we will denote by $\bar{\nu}$:
$${\bar \nu} (
\sum_{\alpha\in\R} A_\alpha
s^\alpha) = \max \{\alpha\;:\; A_{\alpha} \neq  0\}\;.
$$

\subsubsection{Seidel action} \label{subsect-seidel-homom}

Now we are ready to define the Seidel action which is
given by a
homomorphism $\Phi$ from the group $\pi_0 (\cG) = \pi_1(G)$ to the
group of invertible elements of $\oQ$ (see \cite{Se, Lal-McD-Pol})
It sends the class of a loop $\gamma$ to the element

\[
\Phi_\gamma = \sum_{\sigma\in \cS_\gamma} A_\sigma s^{- W
(\sigma)/\Omega}\;.
\]

\subsection{A formula for the restriction of $\tmu$ on
$\pi_1 (G)$}
\label{subsect-explicit-formula-tmu}

\begin{prop}[cf. \cite{Oh-new-1}]
\label{prop-c-on-loops}
Let $[\gamma] \in \pi_1 (G)\subset \tG$ be represented by a loop $\gamma$.
Then for any $a\in QH_{ev} (M)$
\[
c (a, [\gamma]) =
\Omega \bar{\nu} (a\Phi_\gamma^{-1}).
\]

\end{prop}

\bigskip
\noindent
{\bf Proof of Proposition~\ref{prop-c-on-loops}:}
As before let
$K\in \cF$ be the normalized Hamiltonian generating $\gamma$.
Fix a lift of $T_{\gamma}$ associated to some section class
$\sigma$.

Taking $H$ to be the zero Hamiltonian generating the
identity  and
applying formulas (\ref{vsp1}),(\ref{vsp3}) and (\ref{vsp2})
we get that
for any
$a\in QH_{ev} (M)$
$$c (a, [\gamma])  = c (a\Psi^{-1}, \id_\tG) + W(\sigma) =
\Omega\nu (a\Psi^{-1}) +W(\sigma) =$$
$$
 \Omega\bar{\nu} (a\Psi^{-1}
s^{W(\sigma)/\Omega}) = \Omega\bar{\nu} (a\Phi_\gamma^{-1}).
$$
The proposition is proven.
\hfill$\square$
\bigskip

\bigskip
Now we will derive a formula for the restriction of $\tmu$ on
$\pi_1 (G)$.
Recall that the algebra $Q=QH_{ev} (M)$ is assumed to be semi-simple and thus,
as an algebra over $k$, it
decomposes into a direct sum of fields:
\[
Q= Q_1 \oplus ... \oplus Q_d
\]
(see
Section~\ref{3.1}).
Let $e$ be the unit element for the field $Q_1$
involved in the definition of the quasimorphisms $r$ and $\tmu$ (see
Sections~\ref{3.1} and \ref{3.5}) so that
\[
\tmu (\tf) =-\hbox{\rm vol}\, (M)\,\cdot \lim_{m \to
\infty}\frac{c (e, \tf^m )}{m}.
\]
According to Proposition~\ref{prop-c-on-loops},
for any $[\gamma]\in \pi_1 (G)$ one has:
\begin{equation}
\label{eqn-tmu-spectral-numbers-loops} \tmu ([\gamma]) =
-\hbox{\rm vol}\, (M)\,\cdot \Omega\lim_{m \to
+\infty}\frac{\bar{\nu} (e\Phi_\gamma^{-m})}{m}.
\end{equation}

\subsection{Proof of Theorem \ref{1.1.D} in the case $M = \C P^n$}
\label{sect-cp-n-proof}

We consider the manifold $M = \C P^n$ equipped with the
Fubini-Study symplectic form $\omega$.
We need to show that the Calabi quasimorphism $\tmu :\tG\to\R$
descends to a function on $G$, i.e.
$\tmu$ vanishes on $\pi_1 (G)\subset \tG$. Indeed, according to
Proposition~\ref{3.6.A} and
Corollary~\ref{cor-mu-lipschitz},
this would give us a continuous homogeneous
Calabi quasimorphism on $G$.

Example~\ref{2.3.D} tells us that $QH_{ev} (M)$ is a field.
Therefore we can assume that,
in the notation of Section~\ref{subsect-explicit-formula-tmu},
$e=[M]$
and hence, according to
(\ref{eqn-tmu-spectral-numbers-loops})
\[
\tmu ([\gamma]) = -\hbox{\rm vol}\, (M)\, \cdot \Omega \lim_{m\to
+\infty} \frac{\bar{\nu} (\Phi_\gamma^{-m})}{m},
\]
for any $[\gamma]\in\pi_1 (G)$ represented by a loop $\gamma$.
Thus it suffices to prove the following fact:
\begin{equation}
\label{eqn-seidel-claim}
\lim_{m\to +\infty}
\frac{\bar{\nu} (\Phi_\gamma^{-m})}{m}
= 0.
\end{equation}
The proof of (\ref{eqn-seidel-claim}) splits into the following
two propositions.

\begin{prop}
\label{prop-seidel-1}
$\Phi_\gamma$ is a monomial of the type
$\delta A^m s^\beta$ for some
$0\neq\delta\in\C, m\in\Z, \beta\in\R$, where $A\in H_{2n-2} (\C P^n)$
is the hyperplane class.
\end{prop}

\begin{prop}
\label{prop-seidel-2}
If
$\Phi_\gamma = \delta s^\alpha$
for some $0\neq\delta\in\C$
then $\alpha = 0$.
\end{prop}

\bigskip
Postponing the proofs of the propositions we first finish the
proof of
(\ref{eqn-seidel-claim}). Indeed, recall from
Section~\ref{subsect-seidel-homom} that the map
$\gamma\mapsto\Phi_\gamma$
is a homomorphism from $\pi_1 (G)$ to the group of
invertible elements of $\oQ$.
Thus
\[
\Phi_{\gamma^{i}} = \Phi_\gamma^{i}
\]
for any $i\in\Z$. Now using the explicit form of $\Phi_\gamma$
given by Proposition~\ref{prop-seidel-1} and the equality
$A^{-(n+1)} = s [M]$ (see Example~\ref{2.3.D}), we can write:
\begin{eqnarray}
\Phi_{\gamma^{-(n+1)}} = \Phi_\gamma^{-(n+1)} =
\delta^{-(n+1)} A^{-(n+1) m}
s^{-(n+1)\beta} =
\nonumber \\
= \delta^{-(n+1)} [M] s^{-(n+1)\beta + m} = \delta^{-(n+1)}
[M], \nonumber
\end{eqnarray}
where the last equality holds because of
Proposition~\ref{prop-seidel-2}.
But
\[
\bar{\nu} (\delta^{-(n+1)} [M]) = 0.
\]
Hence $\bar{\nu} (\Phi_{\gamma^{-(n+1)}}) = 0$.
Since this is true for every $\gamma$ we get that
$\bar{\nu} (\Phi_{\gamma}^k) = 0$ provided $(n+1)$
divides $k$.
This immediately proves
(\ref{eqn-seidel-claim}).
The proof of  Theorem \ref{1.1.D} in
the case $M = \C P^n$ is finished.
\hfill$\square$

\bigskip
\noindent
{\bf Proof of Proposition~\ref{prop-seidel-1}:}
One needs to show  that there exists at most one class
$\sigma\in {\cS}_\gamma$
such that $GW_\sigma ([M], [M], C) \neq 0$ for some $C\in H_\ast
(M)$. The proof follows from the dimension count.

Indeed, assume $\sigma_1 = \sigma_2 + mS\in {\cS}_\gamma$,
$m\neq 0$, and
$GW_{\sigma_1} ([M], [M], C_1) \neq 0$,
$GW_{\sigma_2} ([M], [M], C_2) \neq 0$ for
some
$C_1, C_2\in H_\ast
(M)$. Since the Gromov-Witten invariants are non-zero the virtual
dimension of the corresponding moduli spaces has to be zero.
Using the formula for the virtual dimension \cite{Se} we
get:
\[
\deg\,
(C_i) = 2n + 4 - 2 c_1^{vert} (\sigma_i ),
\ i=1,2.
\]
Hence, according to (\ref{eqn-u-c-1-sigma-S}),
\[
| \deg\, (C_1) - \deg\, (C_2) | =
2mN = 2m (n+1),
\]
because the minimal Chern number $N$ of $\C P^n$ is $n+1$.
But on the other hand,
$| \deg\, (C_1) - \deg\, (C_2) |$
cannot be bigger than $2n$ which leads us to contradiction unless
$m=0$. The proposition is proven.
\hfill$\square$

\bigskip
\noindent
{\bf Proof of Proposition~\ref{prop-seidel-2}:}
Suppose $\Phi_\gamma = \delta [M] s^\alpha$, i.e. $A_\sigma =
\delta [M]$, $\delta\in\C$. This means that
$A_\sigma\circ_M P = GW_\sigma ([M], [M], P)\neq 0$, where $P
=[point]$.
Let $\cM$ be the moduli space of
$(j,\hJ){\hbox{\rm -holomorphic}}$
sections of the fibration $\pi:E_\gamma\to S^2$ belonging
to the class $\sigma\in {\cS}_\gamma$. We need the following
lemma.

\begin{lemma}
\label{lem-cM-smooth}
$\cM$ is a smooth compact manifold.
\end{lemma}

\smallskip
\noindent {\bf Proof of Lemma~\ref{lem-cM-smooth}:} According to
\cite{Se}, $\cM$ is a smooth manifold of dimension $2n$. The
Gromov compactness theorem \cite{Gro-pshc} says that the only way
the compactness of $\cM$ may fail is a so-called bubbling-off,
when a sequence of sections from $\cM$ converges to a curve in
$E_\gamma$ which is a connected union of a pseudo-holomorphic
section of $\pi$ representing some $\sigma^\prime\in {\cS}_\gamma$
and a number of pseudo-holomorphic spheres lying in fibers of
$\pi$. In such a case the total energy has to be preserved,
meaning that $\sigma = \sigma^\prime + mS$, $m\geq 1$. But, just
as we already checked in the proof of
Proposition~\ref{prop-seidel-1}, the virtual dimension of the
moduli space of pseudo-holomorphic sections belonging to the class
$\sigma^\prime$ equals $2n - 2m (n+1) < 0$ and therefore such
bubbling-off does not happen. Therefore $\cM$ is compact. The
lemma is proven. \hfill$\square$

\bigskip
Consider the
evaluation map
\[
ev: \cM\times S^2\to E_\gamma, \ (v, q) \mapsto v(q).
\]
Then  $\textrm{dim}\, \cM =
2n$ and the degree of the map $ev$ is non-zero
because \break
$GW_\sigma ([M], [M], P)\neq 0$.
Let $\eta$ be the generator of $H^2 (S^2)$ dual to the fundamental class.
Represent $ev^\ast (W)\in H^2 (\cM)\oplus H^2 (S^2)$
as
$ev^\ast (W) = \theta +R\eta$ for some $\theta \in H^2 (\cM)$, $R\in\R$.
Recall that $W^{n+1} = 0$ and that
$\theta^{n+1} = 0$, $\eta^2 = 0$ for dimensional reasons.
Therefore
\[
0 = ev^\ast (W^{n+1}) = (\theta + R\eta)^{n+1} = (n+1)R\theta^n \eta.
\]
On the other hand, the restriction of the coupling class on any
fiber of $\pi$ is the class of the symplectic form on that fiber.
Therefore the product of the $n{\hbox{\rm -th}}$ power of the
coupling class with $\pi^*\eta$ represents a non-zero multiple of the
fundamental class of $E_\gamma$. The image of this cohomology
class under $ev^\ast$ is non-zero, because the degree of $ev$ is
non-zero. Note also that $ev^*\pi^*\eta = \eta$. Thus
\[
0\neq ev^\ast (W^n\pi^*\eta) = (\theta +R\eta)^n \eta =
(\theta^n + nR\theta^{n-1}\eta )\eta = \theta^n\eta.
\]
Combining it with $(n+1) R \theta^n \eta = 0$ we see that $R =
0$.
But
\[
R = ev^\ast (W) ([S^2]) = W (\sigma) = 0.
\]
Now recall that
$\Phi_\gamma = \delta [M] s^{-W (\sigma)/\Omega}$. Hence
$\Phi_\gamma = \delta [M]$. The proposition is proven.
\hfill$\square$

\section{Calabi quasi\-mor\-phism and
combinatorics of level sets of autonomous Hamiltonians
 on the 2-sphere} \prlabel{A}

\subsection{Morse functions and abelian subgroups of $\Ham(S^2)$}\prlabel{A1}

Let $\omega$ be an area form on the 2-sphere $S^2$ with
total area 1. Fix a Morse function $F$ on $S^2$, and consider
the subspace $\cH_F \subset C^{\infty}(S^2)$ consisting of
all functions $H$ whose Poisson bracket with $F$ vanishes:
$\{H,F\} = 0$. For a smooth function $H$ denote by $\psi_H$ the
time-1-map of the Hamiltonian flow generated by $H$. Note that
the Poisson bracket of every two functions from $\cH_F$ vanishes,
and therefore the subgroup
$$\Gamma_F = \{\psi_H\;|\; H \in \cH_F\} \subset \Ham(S^2)$$
is abelian. Intuitively speaking, $\Gamma_F$ is a maximal torus
in $\Ham (S^2)$.

Write $G=\Ham(S^2)$ and let $\mu$ be a continuous homogeneous
Calabi quasimorphism on $G$.
The purpose of this section is to
calculate the restriction of any such $\mu$
on the subgroup $\Gamma_F$. Note that since $\Gamma_F$ is abelian
the map $\mu: \Gamma_F \to \R$ is a homomorphism. It turns
out that the answer can be given in terms of simple
combinatorial data associated to the Morse function $F$
(see Theorem~\ref{answer} below).

\subsection{A measured tree associated to a Morse function}\prlabel{A2}

Let $F: S^2 \to \R$ be a Morse function. Look at connected
components of non-empty level sets of $F$.
These components split into three different groups:

\medskip
\noindent
{\bf I.} Points of local maximum/minimum of $F$;

\medskip
\noindent
{\bf II.} Immersed closed curves whose self-intersections correspond
to critical points
of index 1 of $F$;

\medskip
\noindent
{\bf III.} Simple closed curves.

\medskip
\noindent
We denote by $\cV_1$ and $\cV_2$ the sets of all components of types I
and II respectively. Put $\cV = \cV_1 \cup \cV_2$. Let us emphasize
that $\cV$ is a finite set.

The set $S^2 \setminus \bigcup_{P \in \cV} P$ is a union of a
finite number of pairwise disjoint open cylinders diffeomorphic
to $S^1 \times \R$. Denote their collection by $\cE$.
Every cylinder $C \in \cE$ is foliated by simple closed curves of
type III.  Denote by $e_C$ the space of leaves of this foliation, which is
naturally homeomorphic to $\R$. Consider the {\it Reeb graph} $T$ associated to
the function $F$ as follows (see \cite{Reeb}, cf. \cite{BF}). Its
vertices $v_P$ are in one-to-one correspondence with elements
$P \in \cV$, and its open edges are $e_C$, $C \in \cE$. We say
that an edge $e_C$ connects vertices $v_P$ and $v_Q$ if
$\partial C \subset P \cup Q$. Note that
vertices $v_P$, $P \in \cV_1$
are free. This means that $v_P$ is adjacent to only one edge. If $P \in \cV_2$
the vertex $v_P$ lies in the interior of $T$. Let us emphasize that
every point $x \in T$ corresponds to a subset of $S^2$ which is denoted
by $\gamma_x$.

We claim that in fact $T$ is a tree. Indeed remove a point $x \in T$
lying on any open edge. Since $S^2 \setminus \gamma_x$ is disconnected
we conclude that $T \setminus \{x\}$ is disconnected as well, so $T$ has
no cycles. The claim follows.

Introduce a probability measure $\varrho$ on the tree $T$ which
is uniquely determined by the following conditions. Given
two points $x,y$ lying on the open edge $e_C$, $C \in \cE$, we
define the measure $\varrho ([x,y])$ of the segment $[x,y]$ as
the area of the subcylinder of $C$ bounded by closed curves
$\gamma_x$ and $\gamma_y$. We also require that
all the vertices have measure
$0$.

By definition, a {\it measured tree} is a
finite tree equipped  with a non-atomic Borel
probability measure
whose restriction on every open edge is homeomorphic
to the Lebesgue measure on an open interval. With this language,
the construction above associates a measured tree $(T,\varrho)$
with any Morse function $F$ on $S^2$.

\subsection{The median of a measured tree}\prlabel{A3}

Let $(T,\varrho)$ be a measured tree. A point $x \in T$ is
called {\it a median} if the measure of each connected component
of $T \setminus \{x\}$ does not exceed $\frac{1}{2}$.

\begin{prop}
\prlabel{A3.1}
Every measured tree has unique median.
\end{prop}

\bigskip
\noindent
{\bf Proof:}

\medskip
\noindent
{\bf I. Uniqueness.} Assume on the contrary that $x$ and $y$
are two distinct medians. Denote by $Y$ the connected component
of $T\setminus \{y\}$ which contains $x$. Consider the collection
$X_1, ... , X_m$ of all connected components of $T \setminus \{x\}$.
Assume without loss of generality that $y \in X_1$. Denote by $\alpha$
the open path connecting $x$ and $y$. Then
$$X_2 \cup ... \cup X_m \cup \alpha \subset Y,$$
so
$$\sum_{i=2}^m \varrho(X_i) + \varrho(\alpha) \leq \frac{1}{2}.$$
But $\varrho(X_1) \leq 1/2$ as well, and
$$\sum_{i=1}^m \varrho(X_i) = 1,$$
so
$$\sum_{i=2}^m \varrho(X_i) \geq \frac{1}{2}.$$
This yields $\varrho(\alpha) = 0$, which contradicts to our
assumption that $x\neq y$. Uniqueness follows.

\medskip
\noindent
{\bf II. Existence.} For a point $x \in T$ denote by
$Z_x$ the set of all connected components of $T \setminus \{x\}$.
Put $$\phi (x) = \max_{X \in Z_x} \varrho (X).$$
We claim that $\phi$ is a lower semicontinuous function on $T$.
Obviously, $\phi$ is continuous at $x$ if either $x$ lies on an open edge
of $T$ or $x$ is a free vertex. Suppose that $x$ is an interior vertex
and $Z_x = \{X_1, ... ,X_m\}$. Then
$$\phi (x) = \max_{1 \leq i \leq m} \varrho (X_i).$$
Take sufficiently small
$\epsilon > 0$  and consider
a neighborhood $U$ of $x$ in $T$ consisting of all points $y \in T$
which belong to the edges adjacent to $x$ and satisfy
$\varrho ([x,y]) < \epsilon$. Assume without loss of generality
that $y \in X_1$. Then
$$\phi(y) = \max(\varrho(X_1) - \varrho([x,y]), \sum_{i=2}^m \varrho (X_i) +
 \varrho ([x,y]))$$
for all $y \in U$. We see that $\phi (y) \geq \phi(x) -\epsilon$,
and the claim on the lower semicontinuity of $\phi$ follows.

Since $T$ is compact, the function $\phi$ attains its minimal value
at some point $x \in T$. Let us check that $\phi(x) \leq 1/2$.
Indeed, suppose on the contrary that $\phi(x) = \frac{1}{2} + \delta$
with $\delta > 0 $.
Let $X$ be the (unique)
connected component of $T \setminus \{x\}$
with $\varrho(X) = \frac{1}{2} + \delta$. Denote by $e$ the open
edge of $X$
adjacent to $x$. Choose any point $y \in e$ with $\varrho ([x,y]) < \delta$.
Note that $X \setminus [x,y]$ is a connected
component of $T \setminus \{y\}$ whose measure
equals
$$\varrho(X) - \varrho([x,y]).$$
Since this number is strictly bigger than $\frac{1}{2}$ we
conclude that it is equal to $\phi (y)$. But then
$\phi (y) < \phi (x)$ which contradicts to the assumption that
$\phi$ attains its minimum at $x$. Therefore $\phi(x) \leq 1/2$,
and hence $x$ is a median. This completes the proof of the proposition.
\hfill $\square$

\subsection{The calculation} \prlabel{A4}
Let $\mu: \Ham(S^2) \to \R$ be any continuous homogeneous Calabi
quasimorphism.
For a Morse function $F$ on $S^2$
consider the subgroup $\Gamma_F$  defined in Section~\ref{A2}.
Recall that $\Gamma_F$ consists of all Hamiltonian diffeomor\-phisms
$\psi_H,\; H \in \cH_F$, where the space $\cH_F$ consists of all
functions whose Poisson bracket with $F$ vanishes.
Below we calculate the homomorphism
$\mu : \Gamma_F \to \R$ in terms of the measured tree $(T,\varrho)$
associated to $F$.

To state our result we start with
the following simple observation. Take any $H \in \cH_F$.
Since $\{H,F\} = 0$
the function $H$ is constant on each connected
component of every level set of $F$. Therefore
$H$ descends to a function
$\bar H$ on the tree $T$.
Denote by $x_0$ the median of $(T,\varrho)$.

\begin{thm}
\prlabel{answer}
Let $\mu : \Ham(S^2) \to \R$ be any
homogeneous continuous Calabi quasimorphism. Then
$$\mu(\psi_H) = \int_{S^2} H\cdot \omega - \bar{H}(x_0)$$
for every function $H \in \cH_F$.
\end{thm}

\bigskip
\noindent
{\bf Proof of Theorem~\ref{answer}:} Take a sequence of functions
$w_i : \R \to \R, \; i\in \N$ such that $w_i(s) \equiv \bar{H}(x_0)$ for
$|s- \bar{H}(x_0)| < \frac{1}{i}$ and $w_i$ converges uniformly to
$w(s)=s$ as $i \to +\infty$. Take any $H \in \cH_F$ and put
$H_i = w_i \circ H$. Then the sequence $H_i$ converges uniformly
to $H$. Since $\mu$ is continuous,
one has
$\lim_{i\to +\infty} \mu (\psi_{H_i}) = \mu (\psi_H)$
and hence
it suffices to show that
$$\mu(\psi_{H_i}) = \int_{S^2} H_i\cdot \omega - \bar{H}_i (x_0)$$
for all sufficiently large  $i$.

Denote by $\gamma_{x_0}$ the level set component of $F$
corresponding to the median ${x_0} \in T$. Note that $S^2
\setminus \gamma_{x_0}$ is the disjoint union of a finite number
of open disks which we denote by $U_1, ... , U_m$. By definition
of the median, the area of each $U_j$ does not exceed
$\frac{1}{2}$. Therefore every open subset whose closure lies in
$U_j$ is displaceable. Consider the function $K = H_i -\bar{H}_i
(x_0)$. Note that $K$ can be decomposed as follows: $$K =
K_1+...+K_m, \;\; \text{where}\;\; \text{supp}(K_j) \subset U_j,\;
j= 1, ... ,m.$$ Since $\mu$ is a Calabi quasimorphism we obtain
$$\mu(\psi_{K_j}) = \int_{S^2} K_j \cdot \omega.$$ Note now that
$$\psi_{H_i} = \psi_{K} = \psi_{K_1}\circ ... \circ \psi_{K_m}.$$
Therefore $$\mu(\psi_{H_i}) = \sum_{j=1}^m \mu(\psi_{K_j}) =
\sum_{j=1}^m \int_{S^2}K_j \cdot \omega = \int_{S^2} H_i \cdot
\omega - \bar{H}_i (x_0),$$ because $\mu$ is homogeneous, all
$\psi_{K_j}$ commute and $\int_{S^2} \omega =1$. The proof is
finished. \hfill $\square$

\subsection{Proof of Corollary~\ref{1.3.EE}}
\label{sect-pf-cor-asympt-length}

First assume that $M$ is any of the manifolds listed in
Theorem~\ref{1.1.D} and $\mu$ is the specific continuous homogeneous Calabi
quasimorphism on $G$ that we constructed in Section~\ref{3}.
Let $F$ be an autonomous Hamiltonian on $M$.
Corollary~\ref{cor-mu-lipschitz} immediately yields that
\begin{equation}
\label{eqn-mu-rho}
\frac{|\mu (\psi_F) |}{ \| F\|_{C^0} } \leq \frac{\lim_{m\to
+\infty} \rho (\id, \psi_F^m)/m}{\| F\|_{C^0} } = \zeta (F).
\end{equation}

Now let $M = S^2$, $\int_M \omega = 1$ and let the rest of
the notation be as above.
Let $\cF_{aut}\subset \cF$
denote the set of all autonomous
Hamiltonians on $S^2$ whose integral over $S^2$ is zero.
Consider the set $\cW\subset \cF_{aut}$ of autonomous
Hamiltonians
$F: S^2\to\R$ which satisfy the following properties:

\smallskip
a) $F$ is a Morse function;

\smallskip
b) $0$ is a regular value of $F$;

\smallskip
c) no connected component of $F^{-1} (0)$ divides
$S^2$ into two parts of area $1/2$.

\begin{lemma}
$\cW$ is an open and dense subset in $\cF_{aut}$ with respect
to the $C^\infty{\hbox{\it -topology}}$.
\end{lemma}

The proof is easy and left to the reader.
Now, according to the lemma,
a generic autonomous Hamiltonian $F\in \cF_{aut}$ belongs to $\cW$.
For any such $F$ the value $\bar{F}(x_0)$ of $\bar{F}$ at the median $x_0$
is non-zero. (When we talk
about $\bar{F}(x_0)$ we view $F$ as an element of
$\cH_F$.)
Therefore,
according to
Theorem~\ref{answer}:
\[
\mu(\psi_F) =  - \bar{F}(x_0) \neq 0.
\]
Thus
$|\mu(\psi_F) | > 0$
and hence,
in view of (\ref{eqn-mu-rho}), $ \zeta (F) > 0$ for a generic $F$. This
finishes the proof of Corollary~\ref{1.3.EE}.
\hfill $\square$

\subsection{Proof of Theorem \ref{1.2.5.A}}\prlabel{bcoh}

1) We think of $S^2$ as of the round sphere
$$\{x_1^2 + x_2^2 + x_3^2 = 1\} \subset \R^3$$
endowed with the symplectic form $\omega$ which equals the
spherical area form divided by $4\pi$ (so the total area equals
1). Put $F(x) = \frac{1}{2} x_3$. Note that $F$ is a Morse
function on the sphere, which is in fact the moment map of the
1-turn rotation around the vertical axes. Therefore the measured
Reeb graph of $F$ can be identified with the image of $F$, that is
with the segment $[-\frac{1}{2};\frac{1}{2}]$, endowed with the
Lebesgue measure. Its median of course is the point $0$.

\bigskip
\noindent 2) We start with the case when $M$ is the
annulus $$\{ (p,q) \subset T^*S^1 \;{\Big |} \; 0<p<0.1\}\;.$$ Fix
$\epsilon \in (0;0.1)$. In view of the discussion above there
exists a symplectic embedding $h_{\epsilon}: M \to S^2$ which
sends each circle $\{p = c\}$ to the level set $\{F = c
-\epsilon\}$.

Write $G_M$ for the group of Hamiltonian diffeomorphisms of $M$
generated by Hamiltonians with compact support. The embedding
$h_{\epsilon}$ induces a monomorphism $\phi_{\epsilon} : G_M \to
\Ham (S^2)$. Let $\mu$ be a Calabi quasimorphism on $\Ham (S^2)$.
Put $\mu_{\epsilon} = \mu \circ \phi_{\epsilon}$. Clearly this is
a quasimorphism on $G_M$. Let ${\bar U} \subset M$ be a closed
disk with smooth boundary, whose interior is denoted by $U$. Then
its image $h_{\epsilon} ({\bar U}) $ is displaceable in $S^2$, and
therefore $\mu_{\epsilon}$ coincides with the Calabi homomorphism
on $G_U$.

To show that any finite collection of $\mu_{\epsilon}$'s is
linearly independent over $\R$ we shall proceed as follows. Take
any compactly supported Hamiltonian of the form $H = H(p)$ on $M$.
Denote by $\psi_H$ the corresponding Hamiltonian diffeomorphism.
It follows from Theorem \ref{answer} that $$\mu_{\epsilon}
(\psi_H) = \Ca_M (\psi_H) - H(\epsilon).$$ This immediately yields
the linear independence.

Finally, after the obvious change of parameter, we can assume that
$\epsilon$ runs over $\R$ instead of $(0;0.1)$. This completes the
proof in the case when $M$ is the annulus.

\bigskip
\noindent 3) Let us turn to the case when $M$ is the disk
$$\{(p,q) \in \R^2 \;{\Big |}\; \pi (p^2 + q^2) < 1\}\;.$$ Fix
$\epsilon \in (\frac{1}{2};1)$. There exists a conformally
symplectic embedding $h_{\epsilon}: M \to S^2$ which sends each
circle $\{\pi (p^2 + q^2) = c\},\; c \in(0;1)$ to the level set
$\{F =\frac{1}{2}- \epsilon \cdot c\}.$ The embedding
$h_{\epsilon}$ induces a monomorphism $\phi_{\epsilon} : G_M \to
\Ham (S^2)$. Let $\mu$ be a Calabi quasimorphism on $\Ham (S^2)$.
Then $\mu_{\epsilon} =\epsilon^{-2}\cdot \mu \circ
\phi_{\epsilon}$ is a Calabi quasimorphism on $G_M$.

To show that any finite collection of $\mu_{\epsilon}$'s is
linearly independent over $\R$ we shall proceed as follows. Take
any compactly supported Hamiltonian of the form $H = H{\Big (}\pi
(p^2+q^2){\Big )}$ on $M$. Denote by $\psi_H\in \Ham (M)$ the
corresponding Hamiltonian diffeomorphism. It follows from Theorem
\ref{answer} that $$\mu_{\epsilon} (\psi_H) = \Ca_M (\psi_H) -
\epsilon^{-1} H(\epsilon^{-1}/2).$$ This immediately yields the
linear independence.

Finally, after the obvious change of parameter, we can assume that
$\epsilon$ runs over $\R$ instead of $(\frac{1}{2};1)$. This
completes the proof in the case when $M$ is the disk.
 \hfill $\square$

\bigskip
\noindent
{\bf Acknowledgments.}

\smallskip
We are grateful to P.~Seidel for explaining to us the results
presented in Section \ref{sect-c-on-loops}.
We are indebted to V.~Berkovich
for showing us the proof of Lemma \ref{3.3.A}. We thank D.~McDuff
and Y.-G.~Oh for numerous useful discussions. Part of this work
was done during the stay of the second named author at the
Institute for Advanced Studies (Princeton) in the Winter 2002. He
thanks the Institute and Y.~Eliashberg, the organizer of the
Symplectic Geometry program, for the kind hospitality.

 \bibliographystyle{alpha}

\begin{thebibliography}{99}







\bibitem{Abrams} Abrams, L.,
{\it
The quantum Euler class and
the quantum cohomology of the Grassmannians,
}
Israel J. Math. {\bf 117} (2000), 335-352.



\bibitem{Ban} Banyaga, A.,
{\it
Sur la structure du groupe des diff{\'e}omorphismes qui
pr{\'e}servent une forme symplectique,
}
Comm. Math. Helv. {\bf 53}:2 (1978), 174-227.






\bibitem{Ba-Ghys} Barge, J., Ghys, E.,
{\it
Cocycles d'Euler et de Maslov,
}
Math. Ann. {\bf 294}:2 (1992), 235-265.

\bibitem{Bav} Bavard, C.,
{\it
Longueur stable des commutateurs,
}
L'Enseign. Math. {\bf 37}:1-2 (1991), 109-150.


\bibitem{Bertram-Euler} Bertram, A.,
{\it
Towards a Schubert calculus for maps from
a Riemann surface to a Grassmannian,
}
Internat. J. Math. {\bf 5}:6 (1994), 811-825.


\bibitem{Bertr} Bertram, A.,
{\it
Quantum Schubert calculus,
}
Adv. Math. {\bf 128}:2 (1997), 289-305.



\bibitem{Bertr-CiF-Fult} Bertram, A., Ciocan-Fontanine, I., Fulton, W.,
{\it
Quantum multiplication of Schur polynomials,
}
J. of Algebra {\bf 219}:2 (1999), 728-746.



\bibitem{BF} Bolsinov, A.V. and Fomenko, A.T.,
{\it Exact topological classification of Hamiltonian flows
on smooth two-dimensional surfaces,}
J. of Math. Sciences {\bf 94}:4 (1999), 1457-1476.


\bibitem{Brooks} Brooks, R.,
{\it Some remarks on bounded cohomology,} in
{\it Riemann surfaces and related topics: Proceedings of
the 1978 Stony Brook Conference (State Univ. New York,
Stony Brook, N.Y., 1978)}, pp. 53--63,
Ann. of Math. Stud., 97,
Princeton Univ. Press,
Princeton, N.J., 1981.


\bibitem{Cal} Calabi, E.,
{\it
On the group of automorphisms of a symplectic manifold,
}
in
{\it
Problems in analysis, 1-26.}
Princeton Univ. Press, 1970.


\bibitem{Co-Ze} Conley, C., Zehnder, E.,
{\it
The Birkhoff-Lewis fixed point theorem and a conjecture of V.I.Arnold,
}
Invent. Math. {\bf 73}:1 (1983), 33-49.

\bibitem{DS} {\it Dynamical systems - II.
Ergodic theory with applications to dynamical
systems and statistical mechanics,
}
Ya. G. Sinai ed.,
Encyclopaedia of Mathematical Sciences, 2,
Springer-Verlag, Berlin, 1989.




\bibitem{Entov} Entov, M.,
{\it
Commutator length of symplectomorphisms,
}
preprint,
{\sl
math.SG/0112012},
2001.





\bibitem{Floer} Floer, A.,
{\it
Symplectic fixed points and holomorphic spheres,}
Comm. Math. Phys. {\bf 120}:4 (1989), 575-611.



\bibitem{Gam} Gambaudo, J.-M.,  A talk at the
{\it Workshop on Asymptotic Topology, Foliations and Dynamical
Systems}, M\"unchen, July 2002.




\bibitem{Gouvea} Gouv\^ea, F.Q.,
{\it
$p\text{-adic}$ numbers. An introduction.
}
Springer-Verlag, Berlin, 1997.

\bibitem{Gro-pshc} Gromov, M.,
{\it
Pseudoholomorphic curves in symplectic manifolds,
}
Invent. Math. {\bf 82}:2 (1985), 307-347.

\bibitem{Ho} Hofer, H.,
{\it
On the topological properties of symplectic maps,
}
Proc. of the Royal Soc. of Edinburgh {\bf 115A}:1-2 (1990), 25-28.



\bibitem{Ho-Sa} Hofer, H., Salamon, D.,
{\it
Floer  homology and Novikov rings,
}
in:
{\it
The Floer Memorial Volume, 483-524,}
Progr. Math., 133, Birkh{\"a}user, 1995.


\bibitem{Lal-McD-metr} Lalonde, F., McDuff, D.,
{\it
The geometry of symplectic energy,
}
Ann. of Math. {\bf 141}:2 (1995), 349-371.

\bibitem{Lal-McD-Pol} Lalonde, F., McDuff, D., Polterovich, L.,
{\it
Topological rigidity of Hamiltonian loops and quantum homology,
}
Invent. Math. {\bf 135}:2 (1999), 369-385.

\bibitem{Liu} Liu, G.,
{\it
Associativity of quantum multiplication,
}
Comm. Math. Phys. {\bf 191}:2 (1998), 265-282.






\bibitem{McD-geom-var} McDuff, D.,
{\it
Geometric variants of the Hofer norm,}
preprint, \break
{\sl math.SG/0103089}, 2001.




\bibitem{McD-Sal-pshc} McDuff, D., Salamon, D.,
{\it
$J{\hbox{\it -holomorphic}}$
curves and quantum cohomology,
}
AMS, Providence, 1994.


\bibitem{McD-Sal-sympl-top} McDuff, D., Salamon, D.,
{\it
Introduction to symplectic topology,
}
Oxford University Press, Oxford, 1995.



\bibitem{Oh1} Oh, Y.-G.,
{\it
Symplectic topology as the geometry of action functional. I.
Relative Floer theory on the cotangent bundle,
}
J. Diff. Geom. {\bf 46}:3 (1997), 499-577.


\bibitem{Oh2} Oh, Y.-G.,
{\it
Symplectic topology as the geometry of action functional. II.
Pants product and cohomological invariants,
}
Comm. Anal. Geom. {\bf 7}:1 (1999), 1-54.


\bibitem{Oh-act} Oh, Y.-G.,
{\it
Chain level Floer theory and Hofer's geometry
of the Hamiltonian diffeomorphism group,
}
preprint,
{\sl
math.SG/0104243},
2001.

\bibitem{Oh-new-1} Oh, Y.-G.,
{\it
Normalization of the Hamiltonian and the action spectrum,
}
preprint,
{\sl
math.SG/0206090},
2002.


\bibitem{Oh-new} Oh, Y.-G.,
{\it
Mini-max theory, spectral
invariants and geometry of the Hamiltonian
diffeomorphism group,
}
preprint,
{\sl
math.SG/0206092},
2002.


\bibitem{Ostr} Ostrover, Y.,
{\it
A Comparison of Hofer's Metrics on
Hamiltonian Diffeomorphisms and Lagrangian Submanifolds,
}
preprint,
{\sl
math.SG/0207070},
2002.


\bibitem{PSS} Piunikhin, S., Salamon, D., Schwarz, M.,
{\it
Symplectic Floer-Donaldson theory and quantum cohomology,
}
in:
{\it
Contact and Symplectic Geometry,
171-200,}
Publ. Newton Inst., 8, Cambridge Univ. Press, 1996.



\bibitem{Pol-metr} Polterovich, L.,
{\it
Symplectic displacement energy for Lagrangian submanifolds,
}
Ergodic Th. and Dynam. Syst. {\bf 13}:2 (1993), 357-367.




\bibitem{Pol-masl} Polterovich, L.,
{\it
Hamiltonian loops and Arnold's principle,
}
{\it Topics in singularity theory, 181-187},
Amer. Math. Soc. Transl. Ser. 2, {\bf 180}, AMS, 1997.


\bibitem{Pol-book} Polterovich, L.,
{\it
The geometry of the group of symplectic diffeo\-mor\-phisms,
}
Birkh{\"a}user, 2001.



\bibitem{Post-new} Postnikov, A.,
{\it
Affine approach to quantum Schubert calculus,
}
preprint,
{\sl
math.CO/0205165},
2002.

\bibitem{Reeb} Reeb, G.,
{\it
Sur les points singuliers d'une forme de
Pfaff compl{\`e}tement int{\'e}grable ou d'une
fonction num{\'e}rique,
}
C. R. Acad. Sci. Paris {\bf 222}, (1946), 847-849.


\bibitem{Ru-Ti} Ruan, Y., Tian, G.,
{\it
A mathematical theory of quantum cohomology,
}
Math. Res. Lett. {\bf 1}:2 (1994), 269-278.

\bibitem{Ru-Ti-1} Ruan, Y., Tian, G.,
{\it
A mathematical theory of quantum cohomology,
}
J. Diff. Geom. {\bf 42}:2 (1995), 259-367.

\bibitem{Sal} Salamon, D.,
{\it
Lectures on Floer homology,
}
in:
{\it
Symplectic geometry and topology (Park City, UT, 1997), 143-229,
}
IAS/Park City Math. Ser., 7, AMS, 1999.




\bibitem{Sch-1} Schwarz, M.,
{\it
On the action spectrum for closed symplectically aspherical manifolds,
}
Pacific J. Math.
{\bf 193}:2 (2000), 419-461.

\bibitem{Sch-talks} Schwarz, M.,
{\it
Conference and seminar talks,
}
1999-2001.


\bibitem{Se} Seidel, P.,
{\it
$\pi_1$
of symplectic automorphisms groups and invertibles in quantum homology rings,
}
Geom. and Funct. Analysis {\bf 7}:6 (1997), 1046-1095.


\bibitem{Sieb-Ti} Siebert, B., Tian, G.,
{\it
On quantum cohomology rings of Fano manifolds and a formula of
Vafa and Intriligator,
}
Asian J. Math. {\bf 1}:4 (1997), 679-695.

\bibitem{Vit} Viterbo, C.,
{\it
Symplectic topology as the geometry of generating functions,}
Math. Ann. {\bf 292}:4 (1992), 685-710.



\bibitem{Wi} Witten, E.,
{\it
Two-dimensional gravity and intersection theory on moduli space,
}
Surveys in Diff. Geom. {\bf 1} (1991), 243-310.


\bibitem{Wit-grasm} Witten, E.,
{\it
The Verlinde algebra and the cohomology of the Grassmannian,
}
in
{\it
Geometry, topology and physics},
Conf. Proc. Lect. Notes Geom. Topology, IV,
International Press (1995), 357-422.



















\end{thebibliography}

\bigskip
\noindent
\begin{tabular}{@{} l @{\ \,} l }
Michael Entov  & Leonid Polterovich \\
Department of Mathematics  & School of Mathematical Sciences \\
Technion -- Israel Institute of Technology  & Tel Aviv University \\
Haifa 32000, Israel  & Tel Aviv 69978, Israel \\
{\it e-mail}: entov@math.technion.ac.il  &
 {\it email}: polterov@post.tau.ac.il\\
\end{tabular}

\end{document}